\documentclass[a4paper,12pt,reqno]{amsart}
\usepackage[T1]{fontenc}
\usepackage[utf8]{inputenc}
\usepackage{anyfontsize}
\usepackage[english]{babel}
\usepackage{stmaryrd}
\usepackage{amsthm,amsfonts,amsmath,amssymb}
\usepackage{verbatim}
\usepackage{siunitx}

\usepackage[%
	left=3.5cm,       
	right=3.5cm,      
	top=3.5cm,        
	bottom=3.5cm,     
	heightrounded,    
	bindingoffset=0mm 
]{geometry}

\usepackage[dvipsnames]{xcolor}
\usepackage{graphicx}

\usepackage{xfrac}
\usepackage{functan}
\usepackage{mathtools}
\usepackage{mathrsfs}

\numberwithin{equation}{section}
\allowdisplaybreaks

\newcommand{\R}{\mathbb{R}}

\DeclareMathOperator{\inr}{inr}
\DeclareMathOperator{\interior}{int}

\DeclareMathOperator{\dist}{dist}
\DeclareMathOperator{\reach}{reach}

\DeclareMathOperator{\esssup}{ess\,sup}

\DeclareMathOperator{\Per}{Per}
\DeclareMathOperator{\argmax}{arg\,max}

\usepackage[textwidth=2.5cm]{todonotes}
\usepackage{enumitem}

\usepackage[
colorlinks=true,
urlcolor=teal,
citecolor=magenta,
linkcolor=teal,
linktocpage,
pdfpagelabels,
bookmarksnumbered,
bookmarksopen,
breaklinks=true]
{hyperref}

\usepackage[capitalize,nameinlink,noabbrev]
{cleveref}

\theoremstyle{plain}
	\newtheorem{theorem}{Theorem}[section]
    \newtheorem{lemma}[theorem]{Lemma}
    \AddToHook{env/lemma/begin}{\crefalias{theorem}{lemma}}
    \crefname{lemma}{Lemma}{Lemmas}
	\newtheorem{proposition}[theorem]{Proposition}
    \AddToHook{env/proposition/begin}{\crefalias{theorem}{proposition}}
    \crefname{proposition}{Proposition}{Propositions}
	\newtheorem{corollary}[theorem]{Corollary}
    \AddToHook{env/corollary/begin}{\crefalias{theorem}{corollary}}
    \crefname{corollary}{Corollary}{Corollaries}

\theoremstyle{definition}
	\newtheorem{definition}[theorem]{Definition}
    \AddToHook{env/definition/begin}{\crefalias{theorem}{definition}}
    \crefname{definition}{Definition}{Definitions}
	\newtheorem{example}[theorem]{Example}
    \AddToHook{env/example/begin}{\crefalias{theorem}{example}}
    \crefname{example}{Example}{Examples}
	\newtheorem{remark}[theorem]{Remark}
    \AddToHook{env/remark/begin}{\crefalias{theorem}{remark}}
    \crefname{remark}{Remark}{Remarks}

\title[Monotonicity and evolution]{Monotonicity and evolution 
\\
of total isoperimetric profiles
\\
under Minkowski regularization}%

\author[Z.~Fattah]{Zakaria Fattah}
\address[Z.~Fattah]{
N\^imes Universit\'e, Laboratoire MIPA, Site des Carmes 7, Place Gabriel P\'eri, 30000 N\^imes, France
}
\email{zakaria.fattah@unimes.fr}

\author[I.~Ftouhi]{Ilias Ftouhi}
\address[I.~Ftouhi]{
{N\^imes Universit\'e, Laboratoire MIPA, Site des Carmes 7, Place Gabriel P\'eri, 30000 N\^imes, France}
}
\email{ilias.ftouhi@unimes.fr}

\author[G.~Saracco]{Giorgio Saracco}
\address[G.~Saracco]{
Dipartimento di Matematica e Informatica, Universit\`a di Ferrara, via Machiavelli 30, 44121 Ferrara (FE), Italy
}
\email{giorgio.saracco@unife.it}



\begin{document}


\begin{abstract}
We study the evolution of the total isoperimetric profile of a planar convex body under volume-preserving Minkowski regularization. Several scalar geometric quantities are known to be monotone along this flow; we prove that the total isoperimetric profile is monotone as well, and characterize the equality cases. We also establish an explicit reparametrization law, by which, for nonnegative times, the profile of the regularized body is determined by that of the initial one, and the order of profiles of bodies of equal area is preserved. As consequences, we obtain the uniform convergence of the profile to the one of the ball of equal area; its joint \(\mathrm C^1\) regularity in the volume and the regularization parameter; the fact that it solves a nonlocal quasilinear Hamilton--Jacobi equation; and sharp inequalities for the profile.
\end{abstract}

\maketitle



\section{Introduction}

Let \(K\subset \R^N\) be a convex body, that is, a compact convex set with non-empty interior. Its inradius is defined by
\[
\inr(K)=\sup\{\,r>0:\exists x\in K\ \text{such that } \dist(x; \partial K) \ge r\,\}.
\]
For \(t\ge -\inr(K)\), the \emph{parallel body} of \(K\) is defined by
\begin{equation}
\label{eq:parallel-body}
K_t =
\begin{cases}
K \ominus |t|B, &\text{if } t<0,\\
K \oplus tB, &\text{if } t\ge 0,
\end{cases}
\end{equation}
where \(B\) denotes the closed unit ball centered at the origin,  \(E \oplus F\), respectively \(E \ominus F\), denotes the Minkowski sum, respectively difference, of \(E\) and \(F\).

For convenience, we shall work with \emph{volume-preserving parallel bodies}, that is,
\begin{equation}
\label{eq:vp-pb}
\widetilde{K}_t = \delta_K(t)\, K_t,
\qquad
\forall t>-\inr(K),
\end{equation}
where
\begin{equation}
\label{eq:renormalization}
\delta_K(t)
=
\left(\frac{|K|}{|K_t|}\right)^{\frac{1}{N}}
\end{equation}
is the volume-preserving prefactor. Here, \(|-|\) denotes the usual \(N\)-dimensional Lebesgue measure.

The family \((K_t)_{t}\) has a natural interpretation as a geometric flow. Denoting by \(d_K\) the signed distance from \(\partial K\), negative in the interior of \(K\), for \(t\ge0\) the outer parallel body \(K_t\) is the zero sublevel set of
\[
u(x,t)=\min_{|y-x|\le t}d_K(y),
\]
which, by the Hopf--Lax formula~\cite{BE84}, is the viscosity solution of the eikonal equation \(\partial_tu+|\nabla u|=0\) with initial datum \(d_K\). In the language of the level-set approach~\cite{CGG91,GL26book}, the boundary \(\partial K_t\) moves along its outer normal with unit speed, while inner parallel bodies are obtained by letting the front propagate inward. The renormalization~\eqref{eq:renormalization} makes this flow volume-preserving, in the spirit of the area- and volume-preserving curvature flows~\cite{Gag86,Hui87}. As \(t\) increases, the sets \(\widetilde K_t\) become progressively rounder; along this volume-preserving regularizing flow \(t\mapsto\widetilde{K}_t\), it is then natural to expect several scalar geometric quantities to improve in a monotone fashion. This is indeed the case for: the inradius~\cite{Mat78}; the Cheeger constant in the planar case~\cite{Fto25}; the perimeter~\cite{DL19}; more generally, the quermassintegrals~\cite{HS22,RS22}. Related monotonicity results for anisotropic parallel bodies have also been established in~\cite{Cra21}.

While the results on quermassintegrals~\cite{RS22} hold in any dimension, the Cheeger constant is treated only in the plane~\cite{Fto25}. Here we work in the plane, \(N=2\): the total profile makes sense in any dimension, but our proofs rely on the structure of planar isoperimetric sets established in~\cite{LS22}, and whether and how the results extend to \(N\ge3\) is not
settled by our methods.

\smallskip
The first purpose of this paper is to show that the monotonicity phenomenon propagates from scalar quantities to the full functional level, namely, to the \emph{total isoperimetric profile}
\[
J(v;\widetilde K_t)=\inf\bigl\{\,\Per(E):E\subset \widetilde K_t,\ |E|=v\,\bigr\},
\qquad v\in[0,|K|].
\]
Notice that the profiles are all defined on the fixed interval \([0,|K|]\), independently of \(t\), thanks to renormalization by \(\delta_K(t)\). Here, \(\Per(-)\) denotes the usual variational perimeter in the whole space \(\R^N\). For convex sets, this agrees with \(\mathcal H^{N-1}\) evaluated on their boundary.

The profile \(J(-;\widetilde K_t)\) keeps track of the least perimeter needed to enclose each prescribed volume inside \(\widetilde K_t\). In fact, the scalar quantities above are all read off the profile: the perimeter is \(\Per(\widetilde K_t)=J(|K|;\widetilde K_t)\), the Cheeger constant is \(\inf_{v}J(v;\widetilde K_t)/v\), and the inradius is determined by the largest volume at which the profile still coincides with that of a ball. The profile therefore carries substantially finer information than any of them, and proving its monotonicity requires controlling the constrained minimization problem at \emph{every} volume
scale.

Isoperimetric profiles have already proved to be an effective tool in the study of geometric flows. Hamilton used isoperimetric estimates for the region enclosed by a curve to rule out slowly forming singularities of the curve-shortening flow~\cite{Ham95}. Andrews and Bryan~\cite{AB11} proved that, under the normalized curve-shortening flow, the isoperimetric profile of the enclosed region remains above that of a suitable model solution if it does so initially. There, the profile is the relative one, in which perimeter is measured inside the region, and the argument rests on a differential inequality combined with a comparison principle. In the present setting, both features are different: the total profile measures perimeter in the whole space, so that minimizers meet \(\partial\widetilde K_t\) tangentially rather than orthogonally; monotonicity is obtained by variational means rather than by a maximum principle; and the profile turns out to satisfy an exact equation.

Our first main result, \cref{thm:monotonicity}, establishes that the complete profile \(t\mapsto J(v;\widetilde{K}_{t})\) is monotonically non-increasing for every fixed volume \(v\in[0,|K|]\) and completely characterizes the equality cases. Our proof bypasses direct ambient perimeter estimates by leveraging the structural duality of planar minimizers proven in~\cite{LS22}. 

This theorem has a few direct consequences. First, the isoperimetric profile along the flow \emph{uniformly} converges to the one of the ball of same measure, see \cref{thm:uniform}. Second, we recover the following existing scalar monotonicity results (complete with equality cases) on: the perimeter of the ambient body~\cite{RS22}, see \cref{prop:rs}; the Cheeger constant~\cite{Fto25}, see \cref{prop:fto}. The latter actually extends the monotonicity to \(p\)-Cheeger constants. We further prove that \(\widetilde{K}_t\) becomes calibrable in finite time, even though this is not a direct consequence of our monotonicity result.

\smallskip
The second purpose of this paper is the study of the evolution of the profile \(J(v; \widetilde{K}_t)\) as a function of the pair \((v, t)\). The second main result is a structural reparametrization law for the profile along the regularization; see \cref{thm:RL} and~\eqref{eq:shifted-RL}.
In particular, this yields an exact description of how the underlying variational problem transforms under the flow: for every \(t\ge0\) and \(v\in[0,|K|]\),
\begin{align*}
J\bigl(\Phi_t(v);\widetilde K_t\bigr)
&
=
\delta_K(t)\bigl(J(v;K)+2\pi t\bigr),
\\
\text{where }
\Phi_t(v)
&
=
\delta_K(t)^2\bigl(v+t\,J(v;K)+\pi t^2\bigr),
\end{align*}
and an analogous statement holds for general time pairs (also negative), on the range of volumes which is admissible. This identity allows us to compute explicitly the derivative of \(J(v;\widetilde K_t)\) with respect to \(t\), see \cref{thm:deriv-gen}.

The reparametrization law has several consequences. First, the flow preserves the order of profiles: if the profile of \(K\) lies above that of a body \(L\) of equal area, then the same holds along the flow, that is, for \(\widetilde{K}_t\) and \(\widetilde{L}_t\), see \cref{thm:comparison}. In other words, for \(t\ge0\) the regularization acts on profiles as an autonomous, order-preserving semiflow: the profile of \(\widetilde K_t\) is determined by that of \(K\) alone, by \cref{thm:RL}; every profile decreases along it, by \cref{thm:monotonicity}; and it converges uniformly to the profile of the ball of the same area, by \cref{thm:uniform}.

Second, the map \((v,t)\mapsto J(v;\widetilde{K}_t)\) is shown to be of class \(\mathrm{C}^1\), see \cref{thm:reg-flow}. 
This, paired with the explicit expression of the partial derivative of \(J(v;\widetilde{K}_t)\) with respect to \(t\), proves that, writing \(I(w,t)=J(w;\widetilde K_t)\) and \(p(t)=\Per(\widetilde K_t)\), such map is a classical solution to the nonlocal quasilinear PDE of Hamilton--Jacobi type
\[
\partial_tI=\delta_K(t)\Bigl(2\pi-\frac{p(t)}{2|K|}\,I
-\partial_wI\,\Bigl(I-\frac{p(t)}{|K|}\,w\Bigr)\Bigr),
\]
see \cref{thm:HJ_PDE} and the corresponding transport law~\eqref{eq:HJ-PDE}. Thus, the eikonal Hamilton--Jacobi equation which moves the bodies in space is mirrored, at the level of their profiles, by a Hamilton--Jacobi equation in the volume variable. This latter is nonlocal, as its coefficients depend on the boundary value \(J(|K|;\widetilde K_t)=\Per(\widetilde K_t)\) of the solution itself, and the usual comparison principles for viscosity solutions do not apply directly, see \cref{rem:HJ-perspectives}.

Third, the monotonicity proved in \cref{thm:monotonicity} imposes a sign on this infinitesimal information. Such sign is governed by a concavity property of the volume of parallel bodies, which allows us to completely characterize the equality cases in the differential inequality, see \cref{thm:differential_ineq}. This allows us to extrapolate a new sharp lower bound to the isoperimetric profile of a fixed convex body, see \cref{prop:lower-bound-iso-profile}.

\smallskip
The paper is organized as follows. In \cref{sec:notation}, we fix notation and recall the background material used throughout. In \cref{sec:monotonicity}, we prove the monotonicity theorem and collect its consequences. In \cref{sec:evolution}, we study the evolution of the profile along the flow: we establish the reparametrization law, derive infinitesimal information from it, and collect its consequences, including the order preservation of profiles and the Hamilton--Jacobi equation. Finally, in \cref{sec:non-convex}, we discuss the extent to which these phenomena fail beyond the convex setting.


\subsection*{Acknowledgments}

Part of this research was carried out while the third author was visiting the first and second authors at N\^{i}mes Universit\'e, with additional support from the OcciMath federation. The authors gratefully acknowledge this support. 
I.F.\ was partially supported by a project PEPS JCJC.
Z.F.\ was supported by the French National Research Agency under the France~2030 program, grant ANR-23-EXES-0005 (project Gardener/MASRA), and by the Occitanie Region.
G.S.\ is a member of INdAM--GNAMPA.


\section{Notation and preliminaries}
\label{sec:notation}


\subsection{Notation}
Throughout the paper, \(K\subset \R^2\) denotes a convex body, that is, a compact convex set with non-empty interior.

We define the \emph{volume-preserving parallel bodies} \(\widetilde{K}_t\) for all \(t>-\inr(K)\) as in~\eqref{eq:vp-pb}, and we denote by \(\delta_K(t)\) the \emph{volume-renormalizing prefactor} defined in~\eqref{eq:renormalization}. As mentioned, we shall also call \(\widetilde{K}_t\) the \emph{regularization of \(K\) with parameter \(t\)} and the map \(t\mapsto\widetilde{K}_t\) the \emph{volume-preserving regularizing Minkowski flow}. For convenience, we also set
\[
r_t = \inr(\widetilde{K}_t),
\]
so that \(r_0=\inr(K)\).

We recall that a convex set \(K\) is a \emph{tangential body} (with respect to \(B\)) if and only if it is homothetic to \(K_t\) for all \(t\in(-r_0,0]\), see~\cite[Lem.~3.1.14]{Sch14book}. For the general definition of tangential bodies, we refer the reader to~\cite[Sect.~2.2 and p.~386]{Sch14book}. We then define
\begin{equation}
\label{eq:tauK}
\tau_K
=
\inf\{\, r> -r_0 : K_r \text{ is not a tangential body}\,\}.
\end{equation}
It is well known that \(\tau_K\le 0\), unless \(K\) is a ball, in which case \(\tau_K=+\infty\). Since we work with volume-preserving parallel bodies, we have
\begin{equation}
\label{eq:constancy-interval}
\widetilde{K}_t=\widetilde{K}_s,
\qquad \forall\, s,t\le\tau_K.
\end{equation}

Finally, we define \(\rho_K\) as the maximal interior rolling radius of \(K\), equivalently the smallest radius of curvature of \(K\), namely
\begin{equation}
\label{eq:radius-curvature}
\rho_K
=
\sup\{\, r\in[0,r_0] : \sigma_K+\sigma_K'' \ge r \,\},
\end{equation}
where \(\sigma_K\) is the support function of \(K\), and the inequality in~\eqref{eq:radius-curvature} is understood in the distributional sense when \(K\) is not of class \(\mathrm{C}^2\). The supremum in~\eqref{eq:radius-curvature} is always attained, and it provides a one-sided upper bound on the curvature of \(\partial K\), in the sense that every point of \(\partial K\) supports an interior ball of radius \(\rho_K\). In particular,
\[
K = K_{-r}\oplus rB,
\qquad \forall\, r\in[0,\rho_K].
\]
Recalling that the ball is the unique tangential body of class at least \(\mathrm{C}^1\), if \(K\) is not a ball, one has \(-\rho_K \ge \tau_K\).

For a measurable set \(E\subset \R^2\), \(|E|\) denotes its Lebesgue measure and \(\Per(E)\) its total De Giorgi perimeter in \(\R^2\), that is, \(\Per(E;\R^2)\). For a planar convex body \(K\), this coincides with \(\mathcal{H}^1(\partial K)\), 
where \(\mathcal{H}^1\) stands for the one-dimensional Hausdorff measure.


\subsection{Preliminaries}
In this subsection we collect several known results that will be used throughout the paper.

First, we recall that the Minkowski addition is associative for \(s,t\ge 0\). This composition law implies an analog one on the volume-preserving bodies as follows.

\begin{lemma}
[Composition of outer parallel bodies]
\label{lem:outer-composition} 
Let \(K\subset\R^2\) be a convex body and let \(0\le s<t\). Then, letting \(\theta_s(t)=\delta_K(s)(t-s)\),
\[
\widetilde K_t=\widetilde{(\widetilde K_s)}_{\theta_s(t)}
\qquad
\text{and}
\qquad
\delta_{\widetilde K_s}\bigl(\theta_s(t)\bigr)=\frac{\delta_K(t)}{\delta_K(s)}.
\]
\end{lemma}
\begin{proof}
Since \(0\le s<t\), we have \(K_t=(K_s)_{t-s}\). 
Using that \((\lambda E)_x=\lambda E_{x/\lambda}\), we get
\[
(\widetilde K_s)_{\theta_s(t)}
=
(\delta_K(s)K_s)_{\delta_K(s)(t-s)}
=
\delta_K(s)(K_s)_{t-s}
=
\delta_K(s)K_t.
\]
Hence, the two renormalized bodies coincide and the composition law for prefactors easily follows from the above chain of equalities.
\end{proof}

Second, Matheron's inequality, proved in~\cite{Mat78}, states that
\begin{equation}
\label{eq:matheron}
r_t \ge r_s,
\qquad \forall\, t>s>-r_0.
\end{equation}
Moreover, equality holds if and only if \(t\le \tau_K\).

Third, we recall Steiner's formulas which apply for all \(t\ge -\rho_K\):
\begin{equation}
\label{eq:steiner-area}
|K_t| = |K| + t\Per(K) + \pi t^2,
\end{equation}
and
\begin{equation}
\label{eq:steiner-per}
\Per(K_t) = \Per(K) + 2\pi t.
\end{equation}
See~\cite[Sect.~4.2]{Sch14book} for \(t\ge 0\), and~\cite[Note~6 of Sect.~4.2]{Sch14book} for \(t<0\).

Below \(-\rho_K\) Steiner's formulas break, but the map \(t\mapsto |K_t|\) remains of class \(\mathrm{C}^1\) and its derivative is \(\Per(K_t)\), as proved in the next lemma, where we also collect the regularity properties of the renormalizing prefactor \(\delta_K(t)\).

\begin{lemma}[Regularity properties]
\label{lem:area-C1}
Let \(K\subset\R^2\) be a convex body. Then, the map \(t\mapsto|K_t|\) is of class \(\mathrm{C}^1\) on \((-r_0,+\infty)\) and
\begin{equation}
\label{eq:area-derivative}
\frac{\mathrm d}{\mathrm dt}|K_t| = \Per(K_t).
\end{equation}
Further, \(t\mapsto\delta_K(t)\) is of class \(\mathrm{C}^1\) on \((-r_0,+\infty)\) and
\begin{equation}
\label{eq:delta-derivative}
\delta_K'(t)
=
-\frac{\delta_K(t)^2}{2|K|}\Per(\widetilde{K}_t).
\end{equation}
\end{lemma}

\begin{proof}
For \(t\ge 0\) the claim~\eqref{eq:area-derivative} follows from Steiner's formulas~\eqref{eq:steiner-area}--\eqref{eq:steiner-per}. Let now \(t<0\). Recall that \(\dist(y;\partial K)\) is \(1\)-Lipschitz and its gradient has unit norm almost everywhere. Since \(\{y\in K\,:\,\dist(y;\partial K)\ge |t|\}=K_{t}\) for every \(t\in(-r_0,0]\), the coarea formula gives
\[
|K|-|K_{t}|
=
\int_{t}^{0}\Per(K_{s})\,\mathrm{d}s,
\qquad
\forall\, t\in(-r_0,0].
\]
The map \(t\mapsto K_t\) is continuous with respect to the Hausdorff metric and the perimeter is continuous on convex bodies; hence,~\eqref{eq:area-derivative} follows.

Finally, from~\eqref{eq:renormalization} we get \(2\log\delta_K(t)=\log|K|-\log|K_t|\), whence
\[
\frac{\delta_K'(t)}{\delta_K(t)}
=
-\frac{\Per(K_t)}{2|K_t|}.
\]
The claim follows from the rescaling properties of perimeter and area.
\end{proof}

Fourth, we recall that the total isoperimetric profile of a convex body \(K\subset \R^2\), defined by
\[
J(v;K)=\inf\bigl\{\,\Per(E): E\subset K,\ |E|=v\,\bigr\},
\qquad v\in[0,|K|],
\]
is closely related to the functional
\[
G(\kappa;K)
=
\sup\{\, \kappa |E|-\Per(E): E\subset K,\ |E|\ge \pi r_0^2 \,\},
\qquad \kappa\ge r_0^{-1},
\]
as shown in the next results.
We start by recalling the structure of minimizers of \(G(\kappa;K)\), as proved in~\cite[Thm.~2.3]{LS22}.

\begin{theorem}
\label{thm:evaluation-G}
Let \(K\subset\mathbb{R}^2\) be a convex body. Then:
\begin{enumerate}[label=(\roman*)]
\item if \(\kappa = 1/r_0\), any set attaining \(G(\kappa;K)\) is the convex envelope of two inballs of \(K\);
\item if \(\kappa > 1/r_0\), the unique set attaining \(G(\kappa;K)\) is
\[
K_{-1/\kappa}\oplus \frac{1}{\kappa}B.
\]
\end{enumerate}
In particular,
\begin{equation}
\label{eq:evaluation-G}
G(\kappa;K)=\kappa|K_{-1/\kappa}|-\frac{\pi}{\kappa}.
\end{equation}
\end{theorem}

Given this result, we now state the corresponding information on \(J(v;K)\), obtained by exploiting the fact that \(J(-;K)\) is dual to \(G(-;K)\) in the sense of the Legendre transform on a suitable range; see~\cite[Prop.~4.5]{LS22}.

\begin{theorem}
\label{thm:LS22}
Let \(K\subset\mathbb{R}^2\) be a convex body. The following hold:
\begin{enumerate}[label=(\roman*)]
\item
\label{item:trivial-branch}
if \(v\le \pi r_0^2\), then any set attaining \(J(v;K)\) is a ball of volume \(v\), and
\[
J(v;K)=2\sqrt{\pi v};
\]
\item
\label{item:duality}
if \(v\ge \pi r_0^2\), then
\begin{equation}
\label{eq:duality}
J(v;K)=\sup\{\, \kappa v-G(\kappa;K): \kappa\ge r_0^{-1}\,\};
\end{equation}
\item
\label{item:first-nontrivial-branch}
if \(\pi r_0^2< v\le \pi r_0^2 + 2r_0\mathcal{H}^1(K_{-r_0})\), then any isoperimetric set of volume \(v\) is the convex envelope of two inballs of \(K\);
\item
\label{item:second-nontrivial-branch}
if \(v> \pi r_0^2 + 2r_0\mathcal{H}^1(K_{-r_0})\), there exists a unique \(r=r(v)\in[\rho_K,r_0)\) such that the unique isoperimetric set of volume \(v\) is
\[
E_v=K_{-r}\oplus rB,
\]
that is, the union of all balls of radius \(r=r(v)\) contained in \(K\);
\item
\label{item:convexity}
the function \(J(-;K)\) is concave on \([0,\pi r_0^2]\) and convex on \([\pi r_0^2,|K|]\). In particular, it is affine in the range~\ref{item:first-nontrivial-branch}; 
\item
\label{item:differentiability}
the function \(J(-;K)\) is \(\mathrm{C}^0\) on \([0,|K|]\) and \(\mathrm{C}^1\) on \((0,|K|)\). Further, if \(E_v\) attains \(J(v;K)\) with \(v\neq0,|K|\), then the curvature of \(\partial E_v\cap \interior{(K)}\) equals the derivative of \(J(-;K)\) at \(v\);
\item
\label{item:differentiability-lv}
in particular, in the volume branch in~\ref{item:first-nontrivial-branch}, the derivative equals \(1/r_0\), whereas in the branch in~\ref{item:second-nontrivial-branch}, the derivative equals \(1/r\), where \(r=r(v)\) is the unique radius prescribed there.
\end{enumerate}
\end{theorem}

The structure of isoperimetric sets in the planar convex case described in~\ref{item:trivial-branch}, \ref{item:first-nontrivial-branch}, and \ref{item:second-nontrivial-branch} was first observed in~\cite[Thm.~3.32]{SZ97}. This structure was extended to a larger class of planar sets in~\cite{LS22}, where the regularity properties of \(J(-;K)\) were also investigated. In particular, we refer the reader to~\cite[Thms.~2.3 and~2.4]{LS22} for the shape of the minimizers in the nontrivial branches of~\ref{item:first-nontrivial-branch} and~\ref{item:second-nontrivial-branch}; to~\cite[Prop.~4.5]{LS22} for the duality formula~\eqref{eq:duality}; to~\cite[Thm.~2.5]{LS22} for the convexity and concavity properties in~\ref{item:convexity}; and to~\cite[Rem.~4.6]{LS22} for the differentiability properties in~\ref{item:differentiability} (continuity of the derivative is not explicitly mentioned, but it easily follows from points~\ref{item:trivial-branch}, \ref{item:first-nontrivial-branch}--\ref{item:second-nontrivial-branch} and the explicit knowledge of its value~\ref{item:differentiability}). Point~\ref{item:differentiability-lv} follows by combining the structure in~\ref{item:first-nontrivial-branch} or~\ref{item:second-nontrivial-branch} with~\ref{item:differentiability}.

It is convenient to isolate the radius appearing in \cref{thm:LS22}. Given a convex body \(L\subset\R^2\) and \(v\in(0,|L|)\), we set 
\begin{equation}
\label{eq:RLv}
r_L(v)
=
\bigl(\partial_v J(v;L)\bigr)^{-1},
\end{equation}
and we extend it (continuously) setting \(r_L(0)=0\) and \(r_L(|L|)=\rho_L\).
By \cref{thm:LS22}~\ref{item:trivial-branch} and~\ref{item:differentiability-lv}, the quantity \(r_L(v)\) is precisely the radius of the balls generating the isoperimetric set of volume \(v\in(0,|L|)\), in each of the three branches. In particular, these isoperimetric sets are always a suitable union of balls of radius \(r_L(v)\), hence their interior rolling radius is at least \(r_L(v)\). 

\begin{remark}
\label{rem:properties-rK}
Notice that \(r_L(-)\) is: continuous on \((0,|L|)\) by~\cref{thm:LS22}~\ref{item:differentiability}; nondecreasing on \([0,\pi \inr(L)^2]\); constant in the range of~\ref{item:first-nontrivial-branch}; and nonincreasing in the range of~\ref{item:second-nontrivial-branch}, with \(r_L(|L|)=\rho_L\).
\end{remark}

Directly from \cref{thm:LS22}~\ref{item:convexity}, one can easily prove the following upper bound on the isoperimetric profile; we remark that it is new.

\begin{corollary}
Let \(K\subset\mathbb R^2\) be a convex body different from the ball and let \(v\ge \pi\inr(K)^2\). Then,
\begin{equation}
\label{eq:chord-convexity}
J(v;K) \le \Per(K) - \frac{\Per(K)-2\pi \inr(K)}{|K|-\pi \inr(K)^2}(|K|-v).
\end{equation}
Further, equality for all \(v\ge \pi \inr(K)^2\) is attained uniquely by stadiums.
\end{corollary}

\begin{proof}
The bound~\eqref{eq:chord-convexity} simply follows from convexity of \(J(v;K)\) on the interval \([\pi\inr(K)^2, |K|]\), since the RHS describes the chord between the endpoints of the graph of \(J(-;K)\) between the extrema of said interval.

Equality is achieved if and only if \(J(v;K)\) is linear on the whole interval \(v\ge \pi \inr(K)^2\). This means, that the shape of isoperimetric sets in the nontrivial branch of \(J(v;K)\) is entirely governed by \cref{thm:LS22}~\ref{item:first-nontrivial-branch}, thus \(K\) is the convex envelope of two balls of radius \(\inr(K)\), that is, it is a stadium.
\end{proof}

Also an upper bound only in terms of the perimeter and area of \(K\), that is, independent of the inradius can be proved; also this result appears to be new.

\begin{proposition} 
Let \(K\subset\mathbb R^2\) be a convex body and let \(v\in[0,|K|]\). Then,
\[
J(v;K)\le
\begin{cases}
2\sqrt{\pi v},
&\text{if } v\le
v^*,\\[8pt]
\Per(K)
-\dfrac{\Per(K)+\sqrt{\Per(K)^2-4\pi|K|}}
{2|K|}\,(|K|-v),
&\text{otherwise,}
\end{cases}
\]
where
\[
v^*
=
\dfrac{4\pi |K|^2}
{\bigl(\Per(K)+\sqrt{\Per(K)^2-4\pi|K|}\bigr)^2}.
\]
Equality on the whole range \([0,|K|]\) is uniquely attained by stadiums.
\end{proposition}

\begin{proof}
For the sake of brevity, let \(A=|K|\), \(P=\Per(K)\), and \(r_0=\inr(K)\). Further, let also
\[
r_*=
\frac{P-\sqrt{P^2-4\pi A}}{2\pi},
\]
which is well-defined thanks to the isoperimetric inequality. Notice that
\begin{equation}
\label{eq:r*-v*}
v^* = \pi r_*^2.
\end{equation}

If \(K\) is a ball, \(v^*=|K|\) and the claim is trivial. We may therefore assume that \(K\) is not a ball.

By~\cite[Prop.~2.3]{AMPS22}, we have 
\(\pi r_0^2-Pr_0+A\le0\) and equality is attained only by stadiums. Since \(r_*\) is the smallest root of the equation \(\pi r^2-Pr+A=0\), we have that \(r_*\le r_0\). Hence, the bound holds true for \(v\le v^*=\pi r_*^2\), since here \(J(v;K) = 2\sqrt{\pi v}\).

We now show that for \(v>v^*=\pi r_*^2\), one has
\begin{equation}
\label{eq:intermediate-est-J}
J(v;K) \le P - \frac{1}{r_*}(A-v),
\end{equation}
from which the claim immediately follows from~\eqref{eq:r*-v*} and the definition of \(v^*\).

Let \(v\in(\pi r_*^2, \pi r_0^2]\). By \cref{thm:LS22}, the isoperimetric profile satisfies \(J(v;K)=2\sqrt{\pi v}\), which is concave. Thus, its graph lies below its tangent line at \(v=\pi r_*^2\), so that
\[
J(v;K) \le 2\pi r_* + \frac{v-\pi r_*^2}{r_*} = \pi r_* + \frac{v}{r_*}.
\]
Using that \(r_*\) is a root of \(\pi r^2-Pr+A=0\), inequality~\eqref{eq:intermediate-est-J} follows.

Let now \(v>\pi r_0^2\). We have
\begin{equation}
\label{eq:rPA}
r_*(P-2\pi r_0) - (A-\pi r_0^2)
=
(r_*P -\pi r_*^2 - A) + \pi (r_0-r_*)^2
=
\pi(r_0-r_*)^2
\ge 
0,
\end{equation}
which allows us to estimate from above the negative coefficient of the chord appearing in the RHS of~\eqref{eq:chord-convexity} with \(-1/r_*\), so that~\eqref{eq:intermediate-est-J} follows.

Let us now turn our attention to cases of equality. Independently of \(K\), the bound is attained for \(v\le v^*=\pi r_*^2\). Let then \(v>v^*=\pi r_*^2\).

Since equality in~\eqref{eq:chord-convexity} uniquely holds for stadiums, the bound can be globally attained at most by stadiums. For these sets, one has \(r_*=r_0\). A direct computation exploiting \cref{thm:LS22}~\ref{item:first-nontrivial-branch} shows that they satisfy equality in~\eqref{eq:intermediate-est-J} for all \(v>\pi r_*^2\), whereas equality in~\eqref{eq:rPA} is attained by~\cite[Prop.~2.3]{AMPS22}, so that the claim follows.
\end{proof}

A similar, in spirit, lower bound is proved in \cref{prop:lower-bound-iso-profile}. Its proof relies on the content of \cref{sec:monotonicity,sec:RL}.

\smallskip
As mentioned in the introduction, for ease of notation we state and prove our monotonicity results for the profile of the volume-preserving parallel bodies \(\widetilde{K}_t\). Nevertheless, the following scaling law for the isoperimetric profile allows the reader to translate the results to the profile of the non-rescaled parallel bodies \(K_t\).

\begin{lemma}[Scaling of the profile]
\label{lem:scaling-profile}
Let \(K\subset \R^2\) be a convex body and let \(\alpha>0\). Then
\[
J(\alpha^2 v;\alpha K)=\alpha\, J(v;K),
\qquad \forall\, v\in[0,|K|].
\]
\end{lemma}

\begin{proof}
Let \(\alpha>0\). We have \(E\subset K\) if and only if \(\alpha E\subset \alpha K\). Moreover, by the scaling properties of Lebesgue measure and perimeter,
\[
|\alpha E|=\alpha^2|E|
\qquad\text{and}\qquad
\Per(\alpha E)=\alpha\,\Per(E).
\]
Thus, for any \(F\subset \alpha K\), writing \(F=\alpha E\) for some \(E\subset K\), we obtain
\begin{align*}
J(\alpha^2v;\alpha K)
&
=
\inf\bigl\{\,\Per(F): F\subset \alpha K,\ |F|=\alpha^2v\,\bigr\} 
\\
&
=
\inf\bigl\{\,\Per(\alpha E): E\subset K,\ |E|=v\,\bigr\} 
\\
&
=
\alpha \inf\bigl\{\,\Per(E): E\subset K,\ |E|=v\,\bigr\} 
=
\alpha J(v;K).
\qedhere
\end{align*}
\end{proof}

With this scaling, the monotonicity \cref{thm:monotonicity} can be restated as: for every fixed \(v\in[0,|K|]\), the map
\[
t\mapsto \delta_K(t) J\!\left(\frac{v}{\delta_K(t)^2};K_t\right)
\]
is non-increasing. The constancy part is addressed directly in \cref{thm:monotonicity}.


\section{Monotonicity of the total profile}
\label{sec:monotonicity}

In this section we prove the monotonicity of \(t\mapsto J(v;\widetilde K_t)\) for any fixed \(v\), and we collect some of its consequences.

\subsection{Proof of the monotonicity theorem}
\label{sec:proof-mononicity}

We start by proving the following lemma which is a key ingredient. It follows easily from~\cite[Eq.~(5)]{Fto25}; see also~\cite[Steps~2 and~4 in the proof of Thm.~1.2]{Fto25}.

\begin{lemma}[Area of inner parallel sets]
\label{lem:inner-areas}
Let \(K\subset\mathbb{R}^2\) be a convex body. Then, for all \(t>s>-r_0\), we have
\[
|(\widetilde{K}_t)_{-r}| \ge |(\widetilde{K}_s)_{-r}|,
\qquad \forall r\in[0,r_t].
\]
Moreover, 
for \(r\in(0,r_t)\)
the inequality is strict if and only if \(t>\tau_K\).
\end{lemma}
\begin{proof}
If \(t\le \tau_K\), then \(\widetilde{K}_t=\widetilde{K}_s\), and there is nothing to prove. Henceforth, assume that \(t>\tau_K\). In particular, \(K\) is not a ball and \(\tau_K\le 0\). Since \(\widetilde{K}_\omega\) is constant for all \(\omega\le \tau_K\), we may further assume without loss of generality that \(s\ge \tau_K\). 

Since \(t>s\ge\tau_K\), by Matheron's inequality~\eqref{eq:matheron}, \(r_t>r_s\), thus the inequality is trivially strict in the range \((r_s,r_t)\) and saturated at \(r=r_t\). On the other hand, it is trivially an equality if \(r=0\). Let us now focus on the range of values \((0,r_s]\). 

We first consider the case \(ts\ge 0\). Set
\[
C=\widetilde K_t
\qquad\text{and}\qquad
c=\delta_K(t)(t-s).
\]
First, we claim that \(0<c<\inr(C)\). The lower bound is immediate, since \(t>s\). Regarding the upper bound, as \(s>-r_0\), we have
\[
c 
=
\delta_K(t)(t-s)
<
\delta_K(t)(t+r_0)
=
\inr(C)=r_t,
\]
where the last equality easily follows from the scaling properties of the inradius and its linearity with respect to the Minkowski addition with a ball.

Owing to the scaling rule for parallel sets, \((\lambda E)_{-r} = \lambda E_{-r/\lambda}\), we have
\begin{equation}
\label{eq:C-c}
C_{-c}
=
\bigl(\delta_K(t)K_t\bigr)_{-\delta_K(t)(t-s)}
=
\delta_K(t)(K_t)_{s-t}
=
\delta_K(t)K_s,
\end{equation}
where the last equality holds because we are assuming \(ts\ge0\).

Applying~\cite[Eq.~(5)]{Fto25} to \(C\) and \(C_{-c}\), after rescaling the latter to have area \(|C|\), gives 
\begin{equation}
\label{eq:waldo}
|C_{-r}|
>
\left|
\left(
\left(\frac{|C|}{|C_{-c}|}\right)^{1/2}C_{-c}
\right)_{-r}
\right|,
\qquad
r\in(0,r_s].
\end{equation}
Using now~\eqref{eq:C-c} and that \(|\widetilde{K}_t|=|K|\), we have
\[
\left(\frac{|C|}{|C_{-c}|}\right)^{1/2}C_{-c}
=
\left(
\frac{|K|}{\delta_K(t)^2|K_s|}
\right)^{1/2}
\delta_K(t)K_s
=
\delta_K(s)K_s
=
\widetilde K_s.
\]
The claim now follows by plugging this latter equality in~\eqref{eq:waldo} and recalling that we set \(C=\widetilde{K}_t\).

We are left to prove the inequality in the remaining case \(t>0>s\ge\tau_K\). Applying the previous case first to the pair \((t,0)\) and then to the pair \((0,s)\), the claim follows.
\end{proof}


\cref{lem:inner-areas} implies the monotonicity of the functional \(G(-;\widetilde{K}_t)\), which is dual to \(J(-;\widetilde{K}_t)\), with respect to the regularization parameter \(t\).

\begin{lemma}[Monotonicity of \texorpdfstring{\(G(\kappa;\widetilde{K}_t)\)}{G(kappa; Kt)} with respect to \texorpdfstring{\(t\)}{t}]
\label{lem:monotonicity_G}
Let \(K\subset\mathbb{R}^2\) be a convex body. Then, for all \(t>s>-r_0\), we have
\[
G(\kappa;\widetilde{K}_t)\ge G(\kappa;\widetilde{K}_s),
\qquad \forall \kappa\ge \frac1{r_s}.
\]
Moreover, the inequality is strict if and only if \(t>\tau_K\).
\end{lemma}

\begin{proof}
By Matheron's inequality~\eqref{eq:matheron}, both functionals are well defined on the common range \(\kappa\ge r_s^{-1}\). By \cref{thm:evaluation-G}, letting \(r=\kappa^{-1}\), we have
\[
G(\kappa;\widetilde{K}_t)=\frac{|(\widetilde{K}_t)_{-r}|}{r}-\pi r
\qquad\text{and}\qquad
G(\kappa;\widetilde{K}_s)=\frac{|(\widetilde{K}_s)_{-r}|}{r}-\pi r.
\]
The claim now follows from \cref{lem:inner-areas}.
\end{proof}

\begin{theorem}[Monotonicity of the total profile]\label{thm:monotonicity}
Let \(K\subset\mathbb{R}^2\) be a convex body. Then, for \(t> s > -r_0\), and for every \(v\in[0,|K|]\), we have
\[
J(v;\widetilde{K}_t) \le J(v;\widetilde{K}_s).
\]
Moreover, if \(t\le \tau_K\), then equality holds for all \(v\), whereas if \(t>\tau_K\), the inequality is strict if and only if \(v>\pi r_s^2\).
\end{theorem}

\begin{proof}
If \(K\) is a ball, there is nothing to prove, since \(\widetilde{K}_t=K\) for all \(t\) and \(\tau_K=+\infty\). Henceforth, assume that \(K\) is not a ball; fix \(t>s>-r_0\). We must show that
\begin{equation}
\label{eq:claim}
J(v;\widetilde{K}_t)\le J(v;\widetilde{K}_s),
\qquad \forall v\in[0,|K|].
\end{equation}

If \(t\le \tau_K\), the claim is immediate, because \(\widetilde{K}_t=\widetilde{K}_s\). Thus we may assume \(t>\tau_K\). Since \(J(v;\widetilde{K}_s)\) is constant for \(s\in(-r_0,\tau_K]\), we may further assume without loss of generality that \(s\ge\tau_K\). By Matheron's inequality~\eqref{eq:matheron}, we have the strict inequality
\(
r_s < r_t.
\)

We now split the proof into three regimes, according to the volume \(v\).

\smallskip
\underline{\emph{(i)~Small-volume regime.}}
Let \(v\le \pi r_t^2\). In this range, the isoperimetric set in \(\widetilde{K}_t\) is a ball. Therefore,
\[
J(v;\widetilde{K}_t)
=
2\sqrt{\pi v}
\leq
J(v;\widetilde{K}_s),
\]
and~\eqref{eq:claim} follows. Trivially, as \(r_s<r_t\), equality holds if and only if \(v\le \pi r_s^2\).

\smallskip
\underline{\emph{(ii)~Intermediate-volume regime.}}
Consider now the range of volumes
\[
\pi r_t^2
< 
v
\le 
|(\widetilde{K}_t)_{-r_s}\oplus r_sB|,
\]
which is well defined since \(r_s<r_t\).

On the one hand, we recall that
\[
J(\pi r_s^2;\widetilde{K}_s)
=
2\pi r_s.
\]
For \(v\ge \pi r_s^2\), the profile \(J(v;\widetilde{K}_s)\) is convex by \cref{thm:LS22}\ref{item:convexity}, and its derivative at \(v=\pi r_s^2\) equals \(r_s^{-1}\) by \cref{thm:LS22}\ref{item:differentiability-lv}. Hence,
\[
J(v;\widetilde{K}_s)
\ge 
2\pi r_s+\frac{v-\pi r_s^2}{r_s},
\qquad 
\forall 
v\ge \pi r_s^2.
\]
On the other hand, we also have
\[
J(\pi r_t^2;\widetilde{K}_t)
=
2\pi r_t.
\]
For \(v> \pi r_t^2\), the profile \(J(v;\widetilde{K}_t)\) is convex by \cref{thm:LS22}\ref{item:convexity}; its derivative is increasing in the range under consideration by \cref{thm:LS22}\ref{item:convexity}--\ref{item:differentiability}.
In particular, at \(v=\pi r^2_t\), it is \(r_t^{-1}\). At the other extremum of the interval \(v=|(\widetilde{K}_t)_{-r_s}\oplus r_sB|\) it is at most \(r_s^{-1}\). Indeed: if \(v<|\widetilde{K}_t|\), then it is exactly \(r_s^{-1}\) by \cref{thm:LS22}\ref{item:differentiability-lv}; if instead \(v=|\widetilde{K}_t|\), then the (left) derivative of the profile is \(\rho_{\widetilde{K}_t}^{-1} \le r_s^{-1}\),
where the inequality is a consequence of the fact that \(|\widetilde{K}_t|=|(\widetilde{K}_t)_{-r_s}\oplus r_sB|\).

Therefore,
\[
J(v;\widetilde{K}_t)
\le
2\pi r_t+\frac{v-\pi r_t^2}{r_s}
\]
for every \(v\) in the range under consideration. Thus, it is enough to show that
\[
2\pi r_t+\frac{v-\pi r_t^2}{r_s}
<
2\pi r_s+\frac{v-\pi r_s^2}{r_s},
\]
which is equivalent to
\[
-\frac{\pi}{r_s}(r_t-r_s)^2<0.
\]
This latter holds true, since \(r_s<r_t\).

\smallskip
\underline{\emph{(iii)~Large-volume regime.}}
Finally, let \(v> |(\widetilde{K}_t)_{-r_s}\oplus r_sB|\).
The strict monotonicity in this case follows from the duality with \(G\).

Since \(r_s<r_t\), we are in the nontrivial volume branch introduced in \cref{thm:LS22}\ref{item:second-nontrivial-branch}. Hence, given \(v\) in this range, there exists a unique
\(r=r(v,t)< r_s\)
such that
\[
v = |(\widetilde{K}_t)_{-r}\oplus rB|
\qquad\text{and}\qquad
\argmax_{\kappa\ge r_t^{-1}}\{\kappa v-G(\kappa;\widetilde{K}_t)\}=r^{-1}.
\]
Therefore, by \cref{thm:LS22}\ref{item:duality}, by the above identity, by \cref{lem:monotonicity_G}, and by taking the supremum over \(\kappa\ge r_s^{-1}\), we obtain
\begin{alignat*}{3}
J(v; \widetilde{K}_t)
&
=
&
\sup_{\kappa \ge r_t^{-1}}\{\,\kappa v - G(\kappa; \widetilde{K}_t& ) \,\}
\\
&
=
&
r^{-1} v - G(r^{-1};\widetilde{K}_t)
&
<
r^{-1} v - G(r^{-1};\widetilde{K}_s)
\\
&&
&
\le
\sup_{\kappa \ge r_s^{-1}}\{\,\kappa v - G(\kappa;\widetilde{K}_s) \,\}
=
J(v;\widetilde{K}_s).
\end{alignat*}
This completes the proof.
\end{proof}


\subsection{Consequences of the monotonicity}
\label{sec:consequences-monotonicity}


\subsubsection{Uniform convergence to the profile of the ball}
\label{sec:uniform}

In this subsection we show that the profile of the volume-preserving regularization converges uniformly to the profile of the ball with the same area.

The following lemma holds true, and it is immediate from the fact that: (i) support functions of convex bodies are linear with respect to Minkowski addition, that is,
\(
\sigma_{\alpha K + \beta L} = \alpha \sigma_K + \beta \sigma_L,
\)
as follows from~\cite[Thm.~1.7.5]{Sch14book} and the positive \(1\)-homogeneity of support functions; (ii) Hausdorff convergence is equivalent to uniform convergence of support functions~\cite[Lem.~1.8.14]{Sch14book}.

\begin{lemma}[Asymptotic roundness of \texorpdfstring{\(\widetilde{K}_t\)}{Kt}]
\label{lem:hausdorff}
Let \(K\subset \R^2\) be a convex body, let \((\widetilde{K}_t)_t\) be the volume-preserving parallel bodies of \(K\), and let \(B_K\) be the ball of the same measure as \(K\), centered at the origin. Then, up to a translation, \((\widetilde{K}_t)_t\) converges to \(B_K\) in the Hausdorff metric as $t\rightarrow+\infty$.
\end{lemma}

\begin{theorem}[Uniform convergence of the total profile]
\label{thm:uniform}
Let \(K\subset\mathbb{R}^2\) be a convex body, let \((\widetilde{K}_t)_t\) be the volume-preserving regularization of \(K\), and let \(B_K\) be the ball of the same measure as \(K\), centered at the origin. Then,
\[
J(-;\widetilde{K}_t)\underset{t\to +\infty}{\longrightarrow} J(-;B_K),
\qquad \text{uniformly on }[0,|K|].
\]
\end{theorem}

\begin{proof}
We first show that it is enough to prove pointwise convergence.

Indeed, the isoperimetric profile \(J(-;C)\) is continuous for every convex set \(C\), see \cref{thm:LS22}\ref{item:differentiability}. Moreover, by \cref{thm:monotonicity}, the pointwise convergence is monotone in \(t\). Therefore, by Dini's theorem, pointwise convergence in the compact interval \([0,|K|]\) implies uniform convergence.

We now prove pointwise convergence. By \cref{lem:hausdorff}, the sequence \((\widetilde{K}_t)_t\), up to translation, converges to \(B_K\) in the Hausdorff metric.

\smallskip
\underline{\emph{(i) The case \(v<|K|\).}}
By the Hausdorff convergence, we have \(r_t\to \inr(B_K)\).
Hence, there exists \(\bar t=\bar t(v)\gg 1\) such that
\[
r_t\ge \sqrt{\frac{v}{\pi}},
\qquad \forall t\ge \bar t.
\]
For every such \(t\), we have
\[
J(v;\widetilde{K}_t)=J(v;B_K),
\]
since \(\widetilde{K}_t\) contains a ball of area \(v\).

\smallskip
\underline{\emph{(ii) The case \(v=|K|\).}}
Clearly,
\[
\Per(\widetilde{K}_t)=J(|K|;\widetilde{K}_t)\ge J(|K|;B_K)=\Per(B_K),
\]
since \(\widetilde{K}_t\), respectively \(B_K\), is the only competitor for \(J(|K|;\widetilde{K}_t)\), respectively \(J(|K|;B_K)\). The Hausdorff convergence of \((\widetilde{K}_t)_t\) to \(B_K\), together with convexity of the sets $(\widetilde{K}_t)_t$, implies that \(\Per(\widetilde{K}_t)\to \Per(B_K)\).

The two points above prove the pointwise convergence, and hence the uniform convergence by Dini's theorem.
\end{proof}

\begin{remark}
For future use, we remark that case \textit{(i)} of the above proof implies the following: given any convex body \(K\subset \mathbb{R}^2\) and any \(v\in [0, |K|)\), there exists \(\bar{t}=\bar{t}(v)\) such that
\(
J(v; \widetilde{K}_t) 
=
J(v; B_K)
\)
for every \(t\ge \bar{t}.\)
\end{remark}

\begin{proposition}[Rate of convergence]
Let \(K \subset \mathbb{R}^2\) be a convex body.
Then,
\begin{equation}
\label{eq:rate}
\sup_{v}\bigl(J(v;\widetilde K_t)-J(v;B_K)\bigr)
=
O\!\left(\frac1{t^2}\right),
\end{equation}
as \(t\to+\infty\).
\end{proposition}

\begin{proof}
Being interested in the asymptotics as \(t\to+\infty\), we assume \(t>0\). 

On the one hand, by \cref{thm:LS22}~\ref{item:convexity}, \(J(-;\widetilde K_t)\) is convex on \([\pi r_t^2,|K|]\), differentiable, and
\[
\partial_v J(v;\widetilde K_t)=\frac{1}{r_{\widetilde K_t}(v)},
\]
where \(r_{\widetilde K_t}(v)\in[\rho_{\widetilde K_t},r_t]\) is the radius of the balls generating the isoperimetric set of volume \(v\), see \cref{thm:LS22}~\ref{item:first-nontrivial-branch}--\ref{item:second-nontrivial-branch}.
Hence,
\[
\partial_v J(v;\widetilde K_t)\ge \frac1{r_t},
\qquad\text{for every }v\in[\pi r_t^2,|K|].
\]

On the other hand, by \cref{thm:LS22}~\ref{item:trivial-branch}, the profile \(J(-;B_K)\) has as derivative
\[
\partial_v J(v;B_K)=\frac{1}{r_{B_K}(v)},
\]
where \(r_{B_K}(v)\in(0,(|K|/\pi)^{\sfrac12}]\) is the radius of a ball of measure \(v\). Since \(v\mapsto r_{B_K}(v)\) is increasing,
\[
\partial_v J(v;B_K)\le \frac1{r_t}
\qquad\text{for every }v\in[\pi r_t^2,|K|].
\]
Combining the two inequalities, we obtain
\[
\partial_v\Bigl(J(v;\widetilde K_t)-J(v;B_K)\Bigr)\ge 0,
\qquad\text{for every }v\in[\pi r_t^2,|K|].
\]
Thus, since for \(v\le \pi r_t^2\), the isoperimetric profiles of \(\widetilde{K}_t\) and of \(B_K\) agree, we have
\begin{align}
\sup_{v\ge 0}\bigl(J(v;\widetilde K_t)-J(v;B_K)\bigr)
&
=
\sup_{v\ge\pi r_t^2}\bigl(J(v;\widetilde K_t)-J(v;B_K)\bigr)
\nonumber
\\
&
=
J(|K|;\widetilde K_t)-J(|K|;B_K)
\nonumber
\\
&
=
\Per(\widetilde{K}_t) - \Per(B_K).
\label{eq:uniform-convergence-without-monotonicity}
\end{align}

It remains to compute the asymptotics of the RHS.
By Steiner's formulas~\eqref{eq:steiner-area}--\eqref{eq:steiner-per}, and the scaling of the perimeter, we have
\[
\Per(\widetilde K_t)
=
\delta_K(t)\,\Per(K_t)
=
\sqrt{|K|}\,
\frac{\Per(K)+2\pi t}{\sqrt{|K|+t\Per(K)+\pi t^2}},
\]
whose second order Taylor expansion is
\begin{equation}
\label{eq:second-order-coefficient}
\Per(\widetilde K_t)
=
2\sqrt{\pi |K|}
+
\frac{\sqrt{|K|}}{4\pi^{3/2}}
\bigl(\Per(K)^2-4\pi |K|\bigr)\frac1{t^2}
+o\!\left(\frac1{t^2}\right).
\end{equation}
Since \(B_K\) has perimeter \(2\sqrt{\pi|K|}\), the claim follows.
\end{proof}

\begin{remark}
The previous proposition gives an alternative route to the uniform convergence in \cref{thm:uniform}, independent of the monotonicity theorem. Indeed, the uniform convergence follows from~\eqref{eq:uniform-convergence-without-monotonicity} combined with \cref{lem:hausdorff} and the continuity of the perimeter of convex sets under Hausdorff convergence for convex bodies.
\end{remark}

\begin{remark}
The leading coefficient in~\eqref{eq:second-order-coefficient} is proportional to the isoperimetric deficit of \(K\). Therefore, unsurprisingly, the closer \(K\) is to a ball, the faster the convergence in~\eqref{eq:rate}.
\end{remark}


\subsubsection{Monotonicity of some scalar quantities}
\label{sec:monotonicity-scalar}


A direct application of \cref{thm:monotonicity} allows us to show that several scalar quantities are monotonic decreasing along the flow \(t\mapsto \widetilde{K}_t\).

First, from \cref{thm:monotonicity} applied to \(v=|K|\), 
we are able to recover the two-dimensional statement of~\cite[Thm.~1.1]{RS22}.

\begin{proposition}[Theorem~1.1 of~\cite{RS22}]
\label{prop:rs}
Let \(K\subset\R^2\) be a convex body. Then, for \(t>s>-r_0\), we have
\[
\Per(\widetilde{K}_t)\le \Per(\widetilde{K}_s).
\]
The inequality is strict if and only if \(t>\tau_K\).
\end{proposition}
%


As already noticed in~\cite[Sect.~5.2]{Fto25}, combining \cref{prop:rs} with Matheron's inequality~\eqref{eq:matheron} allows us to infer~\cite[Thm.~1.2]{Lar16}. Thus, it also follows as a corollary of our \cref{thm:monotonicity}, limited to dimension \(2\).

Second, we can prove the monotonicity of the \(p\)-Cheeger constants for \(p\ge 1/2\). The case \(p=1\), corresponding to the usual Cheeger constant, was first proved in~\cite[Thm.~1.2]{Fto25}. We remark that the underlying ingredient in their proof is \cref{lem:inner-areas}, which we also use in a nontrivial way in \cref{thm:monotonicity}.

\begin{proposition}
\label{prop:fto}
Let \(K\subset\R^2\) be a convex body. Then, for \(t>s> -r_0\), we have
\[
h_{p}(\widetilde{K}_t) \le h_{p}(\widetilde{K}_s),
\qquad
\qquad
\forall p\ge 1/2,
\]
where \(h_{p}(-)\) denotes the \(p\)-Cheeger constant. The inequality is strict if and only if \(p>1/2\) and \(t>\tau_K\).
\end{proposition}

\begin{proof}
We recall that the \(p\)-Cheeger constant of a Borel set is defined as
\begin{equation}
\label{eq:cheeger}
h_{p}(\Omega)
=
\inf
\left\{\,
\frac{\Per(E)}{{|E|^p}}
\,:\,
E\subset\Omega,
\,
|E|>0
\,\right\}.
\end{equation}
If \(p=1/2\), the constant reduces to the isoperimetric constant, and the only minimizers are balls within \(\Omega\). Further, if \(K\) is a ball (\(\tau_K=+\infty\)), the inequality is trivially an equality for all \(p, t,\) and \(s\). We therefore restrict ourselves to exponents \(p>1/2\) and the case \(\tau_K \le 0\).

It is easy to see that the infimum can be restricted to isoperimetric sets within \(\Omega\), that is,
\begin{equation}
\label{eq:alternative-cheeger}
h_{p}(\Omega)
=
\inf
\left\{\,
\frac{J(v; \Omega)}{v^p}
\,:\,
v \in (0,|\Omega|]
\,\right\}.
\end{equation}

Let \(C_s\) denote a \(p\)-Cheeger set of \(\widetilde{K}_s\). Using  \(v=|C_s|\) in the alternative formulation~\eqref{eq:alternative-cheeger}, the monotonicity theorem, and finally the minimality of \(C_s\) for \(\widetilde{K}_s\), we obtain
\begin{equation*}
h_p(\widetilde{K}_t)
\le
\frac{J(|C_s|; \widetilde{K}_t)}{|C_s|^p}
\le
\frac{J(|C_s|; \widetilde{K}_s)}{|C_s|^p}
=
\frac{\Per(C_s)}{|C_s|^{p}}
=
h_{p}(\widetilde{K}_s).
\end{equation*}

If \(t\le \tau_K\), then \(\widetilde{K}_t=\widetilde{K}_s\), and equality holds. Let us now assume \(t>\tau_K\). %
We claim that the second inequality is strict. First, observe that since \(K\) is not a ball also \(\widetilde{K}_s\) is not. Therefore, by~\cite[Thm.~2.5(3ii)--(3iv)]{CS23},  \(|C_s|>\pi r_s^2\). \cref{thm:monotonicity} allows us to conclude.
\end{proof}

While the following proposition is not a consequence of monotonicity, we believe it is worth mentioning, as it could have possible applications. By definition of Cheeger constant, we have
\begin{equation}
\label{eq:cheeger-est-flow}
h_1(\widetilde K_t)
\le
\frac{\Per(\widetilde K_t)}{|\widetilde{K}_t|},
\qquad
\forall t>-r_0.
\end{equation}
By \cref{lem:hausdorff} combined with the continuity of the Cheeger constant with respect to the Hausdorff metric~\cite[Prop.~3.1]{Par17}, the LHS converges to \(h_1(B_K)\), whereas by the continuity of the perimeter of convex bodies with respect to the same metric, the RHS to the ratio \(\Per(B_K)/|B_K|\). Since \(B_K\) realizes \(h_1(B_K)\), at the limit, equality in~\eqref{eq:cheeger-est-flow} is attained. Therefore, it is reasonable to wonder whether equality in~\eqref{eq:cheeger-est-flow} could already be attained at finite time, that  is, there exists \(T<+\infty\) such that
\[
h_1(\widetilde K_t)
=
\frac{\Per(\widetilde K_t)}{|\widetilde{K}_t|},
\qquad
\forall t \ge T.
\]
Sets satisfying the above equality are either called \emph{self-Cheeger sets} or \emph{calibrable sets} and are of interest in several applications, such as image reconstruction, see~\cite{ACC05,ROF92}.
The following proposition shows this is the case. Further, \(T\) can be chosen as the smallest time such that \(\widetilde{K}_t\) is calibrable.

\begin{proposition}[Finite-time calibrability]
\label{prop:calibrability}
Let \(K\subset\mathbb{R}^2\) be a convex body. Further, let 
\begin{equation}
\label{eq:first-time-calibrability}
T = -\rho_{K} + \sqrt{\frac{|K_{-\rho_{K}}|}{\pi}}.
\end{equation}
Then,  if \(\rho_{K}<r_0\), \(\widetilde{K}_t\) is calibrable for all \(t\ge T\). If otherwise \(\rho_{K}=r_0\), the same holds true with a strict inequality.
\end{proposition}
\begin{proof}
We recall that a two-dimensional convex body \(C\) is calibrable if and only if
\[
\esssup_{\partial C} \kappa_C \le \frac{\Per(C)}{|C|},
\]
where \(\kappa_{C}\) denotes the curvature of \(C\), refer to~\cite[Thm.~2]{KLR06} (see also~\cite[Crit.~1.5]{Sar21} for the more general sharp statement). This condition can be restated in terms of the smallest radius of curvature of \(C\), \(\rho_{C}\ge 0\), as
\begin{equation}
\label{eq:calibrability-radius}
\rho_{C} \ge  \frac{|C|}{\Per(C)}.
\end{equation}
As calibrability is not affected by rescaling, it is enough to verify it holds for \(K_t\). Testing the above on \(K_t\), and owing to the fact that Minkowski addition by a ball of radius \(t\) increases the radius of curvature exactly by \(t\), we find that~\eqref{eq:calibrability-radius} is satisfied if and only if
\[
(\rho_{K} + t)\Per(K_t)
\ge
|K_t|.
\]
Choosing \(t\ge-\rho_K\), using Steiner's formulas~\eqref{eq:steiner-area}--\eqref{eq:steiner-per}, we find that this is equivalent to
\[
\rho_{K} (\Per(K)+2\pi t) +\pi t^2 \ge |K|.
\]
Using once again Steiner's area formula~\eqref{eq:steiner-area} for \(K_{-\rho_K}\), this is equivalent to
\[
\pi(t+\rho_K)^2 \ge |K_{-\rho_K}|.
\]
Rearranging the claim follows.
\end{proof}


\section{The evolution law}
\label{sec:evolution}

In this section, we derive the exact reparametrization law and then differentiate it with respect to time to obtain infinitesimal monotonicity information.

For the sake of notation, we let
\begin{equation}
\label{eq:Deltatheta}
\theta_s(t) =\delta_K(s) (t-s).
\end{equation}
When \(s=0\) and the base body is \(K\), the increment satisfies
\(
\theta_0(t) = t,
\)
and we will simply use \(t\).
For \(t,s>-\rho_K\), the quantity \(\theta_s(t)\) represents the time increment needed to pass from the convex body \(\widetilde{K}_s\) to \(\widetilde{K}_t\), see the discussion in \cref{ssec:degeneracy-rho}; below this threshold it retains this role in all the formulas of this section, although it is no longer the parameter of an actual Minkowski operation.

For \(t>-\inr(L)\) we also define the set of \emph{admissible volumes}
\begin{equation}
\label{eq:admissible}
\mathcal{A}_t(L)
=
\bigl\{\, v\in[0,|L|] \,:\, r_L(v)\ge \max\{0,-t\} \,\bigr\}.
\end{equation}
where \(r_L(\cdot)\) has been defined in~\eqref{eq:RLv}. Recalling the properties of \(r_L(\cdot)\) stated in \cref{rem:properties-rK}, \(\mathcal{A}_t(L)\) is a closed
interval. Explicitly
\begin{equation}
\label{eq:admissible-explicit}
\mathcal{A}_t(L)
=
\begin{cases}
[0,|L|], & \text{if } t\ge 0,
\\[2pt]
\bigl[\,\pi t^2,\ |L_{-|t|}\oplus|t|B|\,\bigr], & \text{if } t<0;
\end{cases}
\end{equation}
in the case \(t<0\), the right extremum is \(|L|\) if and only if \(|t|\le\rho_L\).
\begin{remark}
\label{rem:use-At}
Let \(E_v\) be an isoperimetric set in \(L\) of volume \(v\in\mathcal{A}_t(L)\).
On the one hand, if \(t\ge 0\), every volume is admissible.
On the other hand, if \(t<0\), the condition \(v\in\mathcal A_t(L)\) ensures, owing to \cref{thm:LS22}, that one has \((E_v)_t \neq \emptyset\) and \(E_v=(E_v)_{-|t|}\oplus|t|B\).
In both cases, one can use Steiner's formulas~\eqref{eq:steiner-area}--\eqref{eq:steiner-per} to find the area and perimeter of \((E_v)_{-|t|}\) from that of \(E_v\).
\end{remark}


\subsection{Transport of balls and of isoperimetric sets}
\label{sec:trasports-balls}
We begin with some preliminary results, which describes how balls are transported along the flow, and, consequently, how isoperimetric sets are.

\begin{lemma}[Ball-to-ball lemma]
\label{lem:ball-to-ball}
Let \(K\subset\mathbb{R}^2\) be a convex body and let \(s,t>-r_0\). A ball \(\rho B(x)\), with \(\rho\ge0\) and \(\rho+\theta_s(t)>0\), is contained in \(\widetilde K_s\) if and only if
\begin{equation}
\label{eq:radii-generic}
\frac{\delta_K(t)}{\delta_K(s)}\bigl(\rho+\theta_s(t)\bigr)\, B\!\left(\frac{\delta_K(t)}{\delta_K(s)}\,x\right)
\subset
\widetilde K_t .
\end{equation}
In particular, if \(s<t\), the constraint \(\rho+\theta_s(t)>0\) is satisfied by every \(\rho\ge0\), and the map
\[
\rho\mapsto \frac{\delta_K(t)}{\delta_K(s)}(\rho+\theta_s(t))
\]
sends \([0,r_s]\) onto \([-\theta_t(s),r_t]\).
\end{lemma}

\begin{proof}
We split the proof in two steps; in the first one, we assume one of the parameters to be zero.

\smallskip
\underline{\emph{Step~(i): Base body \(K\).}} We consider the pair \(s=0\) and \(t>-r_0\). Since \(y\mapsto \delta_K(t)y\) is a homothety, it is enough to prove the analog claim for non-rescaled parallel bodies, that is, \(rB(x) \subset K\) if and only if one has \((r+t)B(x) \subset K_t\). 

On the one hand, if \(t\ge 0\), then
\[
rB(x)\subset K
\quad
\implies
\quad
rB(x)\oplus tB
\subset K\oplus tB=K_t.
\]
The converse implication holds true, since for convex sets and \(t\ge 0\) one has \(K_t\ominus tB = K\).
The claim follows since \(rB(x)\oplus tB\) is exactly \((r+t)B(x)\).

On the other hand, if \(t<0\), we have
\[
y\in K_t = K\ominus |t|B
\quad
\Longleftrightarrow
\quad
\{y\}\oplus|t|B \subset K.
\]
By assumption, \(r+t = r-|t| > 0\), so that given any \(x\), the set \((r-|t|)B(x)\) is a well-defined ball. Therefore,
\[
rB(x)
=
(r-|t|)B(x) \oplus |t| B
\subset K
\quad 
\Longleftrightarrow
\quad
(r-|t|)B(x)
\subset K_t.
\]

\smallskip
\underline{\emph{Step~(ii): Generic time pairs.}}
There is no loss of generality in assuming \(s<t\) as for the opposite cases it is enough to reverse the roles of \(s\) and \(t\). In particular \(\theta_s(t)>0\).

\smallskip
\underline{\emph{Subcase~(i): \(s\le0\).}} We apply twice Step~\textit{(i)}. 
A ball of radius \(\rho\) is contained in \(\widetilde K_s\) if and only if the ball of radius \(\rho/\delta_K(s)-s\) is contained in \(K\). This operation maps \((0,r_s]\) onto \((|s|, r_0]\).
Then, since \(|s|=-s>-t\), we can apply Step~(i) to all \(r\in (|s|,r_0]\). 
For \(\rho=0\), equivalently \(r=|s|\), \(x\in\widetilde K_s\) if and only if \(x/\delta_K(s)\in K_s\), that is, if and only if \(|s|B(x/\delta_K(s))\subset K\).  Applying Step~\textit{(i)} at time \(t\), this is equivalent to
\(-\theta_t(s)B\!\left(\delta_K(t)x/\delta_K(s)\right)\subset \widetilde K_t\).

\smallskip
\underline{\emph{Subcase~(ii): \(s>0\).}}
Passing first from \(\widetilde K_s\) back to \(K\) does not cover all radii in \((0,r_s]\), since the corresponding radius in \(K\) would be \(\rho/\delta_K(s)-s\), which is non-positive for sufficiently small \(\rho\). Since \(\theta_s(t)>0\), we instead apply Step~\textit{(i)} directly to \(\widetilde K_s\) with parameter \(\theta_s(t)\). Together with \cref{lem:outer-composition}, this gives the claim. Note that, as \(\theta_s(t)>0\), the constraint \(\rho+\theta_s(t)>0\) is satisfied for all \(\rho\ge0\), so no degenerate case arises here.
\end{proof}

As a consequence of \cref{lem:ball-to-ball}, of the structure of isoperimetric sets provided by \cref{thm:LS22}~\ref{item:trivial-branch},~\ref{item:first-nontrivial-branch}, and~\ref{item:second-nontrivial-branch}, and of the characterization of the derivative of the profile as the curvature given in \cref{thm:LS22}~\ref{item:differentiability}, the following holds.

\begin{corollary}
\label{cor:iso-one-to-one}
Let \(K\subset\mathbb{R}^2\) be a convex body and let \(s,t>-r_0\). Then, the map
\begin{equation}
\label{eq:iso-transformation}
E\longmapsto\frac{\delta_K(t)}{\delta_K(s)}\,E_{\theta_s(t)}
\end{equation}
is a one-to-one correspondence between the isoperimetric sets of
\(\widetilde K_s\) of volume in the interior of \(\mathcal{A}_{\theta_s(t)}(\widetilde K_s)\) and
the isoperimetric sets of \(\widetilde K_t\) of volume in the interior of
\(\mathcal{A}_{\theta_t(s)}(\widetilde K_t)\).

Under this correspondence, the free boundary curvatures \(\kappa_s\) and
\(\kappa_t\) of two corresponding isoperimetric sets are related by
\begin{equation}
\label{eq:curvature-transformation}
\kappa_t
=
\frac{\delta_K(s)}{\delta_K(t)}\cdot\frac{\kappa_s}{1+\theta_s(t)\,\kappa_s}.
\end{equation}
\end{corollary}
\begin{remark}
\label{rem:endpoints}
The restriction to interior volumes is only due to the degenerate behavior at the endpoints of the admissible ranges and it may be removed.

Isoperimetric sets are determined up to sets of measure zero and, by \cref{thm:LS22}, they admit a convex representative. We choose it closed for all volumes and we choose as the isoperimetric set of volume zero any singleton \(\{x\}\) in the body. With this convention,~\eqref{eq:iso-transformation} is meaningful also at the
endpoints.

As for the curvature relation~\eqref{eq:curvature-transformation}, consider two corresponding isoperimetric sets under~\eqref{eq:iso-transformation}. If both free-boundary radii are finite and positive, then~\eqref{eq:curvature-transformation} holds in the usual sense. There are two degenerate cases to be addressed, where we use the convention that \(1/0=+\infty\), in such a way that when \(r_L(v)=0\), the corresponding curvature is \(+\infty\).

The initial radius could vanish: this occurs for \(v=0\), and also for \(v=|K|\) when \(\rho_{\widetilde K_s}=0\). Since \(v\in\mathcal A_{\theta_s(t)}(\widetilde K_s)\), we necessarily have
\(\theta_s(t)\ge0\). For \(\theta_s(t)>0\) the RHS of~\eqref{eq:curvature-transformation} is then understood by letting \(\kappa_s\to+\infty\); the case \(s=t\) is trivial.

Conversely, the transported radius could vanish: this happens if the initial radius is positive but equals \(|\theta_s(t)|\) with \(\theta_s(t)<0\). Hence, the denominator in~\eqref{eq:curvature-transformation} is zero, while the numerator is finite, so that \(\kappa_t=+\infty\).
\end{remark}


\subsection{The reparametrization law}
\label{sec:RL}
Here we prove that the isoperimetric profile of \(\widetilde K_t\) is obtained from that of \(\widetilde{K}_s\) through a suitable reparametrization, for every \(t,s>-r_0\).

Before stating the main theorem of this section, it is convenient to define the map that transports volumes along the flow. 
Given \(s,t>-r_0\), we let
\[
\Phi_{s,t}\colon
\mathcal{A}_{\theta_s(t)}(\widetilde K_s)
\longrightarrow
\mathcal{A}_{\theta_t(s)}(\widetilde K_t)
\]
be defined as
\begin{equation}
\label{eq:Phi}
\Phi_{s,t}(v)
=
\frac{\delta_K(t)^2}{\delta_K(s)^2}
\Bigl(v+\theta_s(t)\,J(v;\widetilde K_s)+\pi\,\theta_s(t)^2\Bigr),
\end{equation}
where the intervals of admissible volumes have been introduced in~\eqref{eq:admissible}. By \cref{cor:iso-one-to-one} and \cref{rem:use-At}, \(\Phi_{s,t}\) sends the volume of an isoperimetric set of \(\widetilde K_s\) to the volume of the corresponding isoperimetric set of \(\widetilde K_t\). We collect below a few 
properties of \(\Phi_{s,t}\) we will use throughout.

\begin{lemma}
[Properties of \texorpdfstring{\(\Phi_{s,t}\)}{Phi}]
\label{lem:Phi}
Let \(K\subset\R^2\) be a convex body and let \(s>-r_0\) be fixed. Set
\begin{equation}
\label{eq:domain-Phi}
\mathcal Q_s
=
\left\{
(v,t)\in[0,|K|]\times(-r_0,+\infty):
v\in\mathcal A_{\theta_s(t)}(\widetilde K_s)
\right\}.
\end{equation}
Then, 
\begin{enumerate}[label=(\roman*)]
\item
\label{item:phi-continuous}
the map \((v,t)\mapsto\Phi_{s,t}(v)\) is continuous on \(\mathcal Q_s\);
\item
\label{item:phi-identity}
\(v\mapsto \Phi_{s,s}(v)\) is the identity;
\item
\label{item:phi-v-derivative}
the partial derivative \(\partial_v \Phi_{s,t}\) satisfies
\begin{equation}
\label{eq:Phi-derivative}
\partial_v\Phi_{s,t}(v)
=
\frac{\delta_K(t)^2}{\delta_K(s)^2}
\Bigl(1+\theta_s(t)\,\partial_vJ(v;\widetilde K_s)\Bigr),
\end{equation}
and it is strictly positive on the interior of
\(\mathcal{A}_{\theta_s(t)}(\widetilde K_s)\);
\item
\label{item:phi-increasing}
for each fixed \(t>-r_0\), the map \(v\mapsto\Phi_{s,t}(v)\) is a strictly increasing homeomorphism of \(\mathcal{A}_{\theta_s(t)}(\widetilde K_s)\) onto \(\mathcal{A}_{\theta_t(s)}(\widetilde K_t)\).
\end{enumerate}
Further,
\begin{enumerate}[label=(\roman*),resume]
\item
\label{item:phi-IFT-domain}
for every fixed \(w_0\in(0,|K|)\), there exists an open neighborhood \(U\) of \((w_0,s)\) in \((0,|K|)\times(-r_0,+\infty)\) such that, for all \((w,t)\in U\), one has \(w\in\mathcal A_{\theta_t(s)}(\widetilde K_t)\) and \(\Phi_{s,t}^{-1}(w)\) lies in the interior of \(\mathcal A_{\theta_s(t)}(\widetilde K_s)\);
\item
\label{item:phi-IFT-C1}
the map \((w,t)\mapsto\Phi_{s,t}^{-1}(w)\) is of class \(\mathrm C^1\) on \(U\).
\end{enumerate}
\end{lemma}

\begin{proof}
Point~\ref{item:phi-identity} is immediate since \(\theta_s(s)=0\).
By \cref{thm:LS22}~\ref{item:differentiability} the map \(v\mapsto J(v;\widetilde K_s)\) is continuous, whereas \(t\mapsto\delta_K(t)\) is continuous by \cref{lem:area-C1}. 
Recalling the definition of \(\theta_s(t)\) given in~\eqref{eq:Deltatheta}, the explicit formula~\eqref{eq:Phi} is continuous; thus its restriction to \(\mathcal Q_s\) is as well, settling~\ref{item:phi-continuous}.

Point~\ref{item:phi-v-derivative} follows by differentiating~\eqref{eq:Phi} with respect to \(v\). Recalling now the definition of \(r_{\widetilde K_s}\), see~\eqref{eq:RLv}, and that of \(\mathcal{A}_{\theta_s(t)}(\widetilde K_s)\), see~\eqref{eq:admissible}, this derivative is strictly positive on the interior of \(\mathcal{A}_{\theta_s(t)}(\widetilde K_s)\). 
Indeed, if \(\theta_s(t)\ge 0\), the claim is immediate. If \(\theta_s(t)<0\), then \(r_{\widetilde K_s}(v)>|\theta_s(t)|\) for every \(v\) in the interior of \(\mathcal A_{\theta_s(t)}(\widetilde K_s)\). Hence,
\[
1+\theta_s(t)\,\partial_v J(v;\widetilde K_s)
=
1-\frac{|\theta_s(t)|}{r_{\widetilde K_s}(v)}
>0.
\]
Thus, \(\Phi_{s,t}\) is continuous and strictly increasing on the closed interval \(\mathcal{A}_{\theta_s(t)}(\widetilde K_s)\), and therefore a homeomorphism onto its image, which is a closed interval. By \cref{cor:iso-one-to-one} and \cref{rem:use-At}, \(\Phi_{s,t}\) maps the interior of \(\mathcal{A}_{\theta_s(t)}(\widetilde K_s)\) onto the interior of \(\mathcal{A}_{\theta_t(s)}(\widetilde K_t)\). Two closed intervals with the same interior coincide, so the image is \(\mathcal{A}_{\theta_t(s)}(\widetilde K_t)\), proving~\ref{item:phi-increasing}.

To prove points~\ref{item:phi-IFT-domain}--\ref{item:phi-IFT-C1}, first notice that \(r_{\widetilde K_s}(w_0)>0\), since \(w_0\in(0,|K|)\). By continuity of \(r_{\widetilde K_s}\) at \(w_0\), and since \(\theta_s(t)\to0\) as \(t\to s\), there exists an open neighborhood of \((w_0,s)\) such that \(r_{\widetilde K_s}(v)>\max\{0,-\theta_s(t)\}\). Hence \(v\) lies in the interior of \(\mathcal A_{\theta_s(t)}(\widetilde K_s)\). Hence, \((v,t)\mapsto\Phi_{s,t}(v)\) is well-defined in such a neighborhood and it is of class \(\mathrm C^1\) as one can see directly from~\eqref{eq:Phi} combined with \cref{lem:area-C1} and \cref{thm:LS22}~\ref{item:differentiability}.

We can thus consider the map \(F(v,w,t)=\Phi_{s,t}(v)-w\). Since \(\Phi_{s,s}\) is the identity, we have \(F(w_0,w_0,s)=0\) and \(\partial_vF(w_0,w_0,s)=1\). By the Implicit Function Theorem, there exist an open neighborhood \(U\) of \((w_0,s)\) and a map \(\Psi\in\mathrm C^1(U)\) such that \(\Psi(w_0,s)=w_0\) and \(\Phi_{s,t}(\Psi(w,t))=w\) for every \((w,t)\in U\).

Up to shrinking \(U\), we may assume that \(\Psi(w,t)\) lies in the interior of \(\mathcal A_{\theta_s(t)}(\widetilde K_s)\). By~\ref{item:phi-increasing}, this implies that \(w=\Phi_{s,t}(\Psi(w,t))\in\mathcal A_{\theta_t(s)}(\widetilde K_t)\). Since \(\Phi_{s,t}\) is injective on \(\mathcal A_{\theta_s(t)}(\widetilde K_s)\), we conclude that \(\Psi(w,t)=\Phi_{s,t}^{-1}(w)\). This proves~\ref{item:phi-IFT-domain}, while~\ref{item:phi-IFT-C1} follows from the fact that \(\Psi\in\mathrm C^1(U)\).
\end{proof}

We are now ready to state and prove the main theorem of this section, namely, the reparametrization law, which relates the isoperimetric profiles at two arbitrary times.
\begin{theorem}[The profile reparametrization law]
\label{thm:RL}
Let \(K\subset\R^2\) be a convex body, let \(s,t>-r_0\), and let
\(\theta_s(t)=\delta_K(s)(t-s)\) be as in~\eqref{eq:Deltatheta}. Then,
\begin{align}
\tag{RL}
\label{eq:shifted-RL}
J(\Phi_{s,t}(v);\widetilde{K}_t)
=
\frac{\delta_K(t)}{\delta_K(s)}\big(J(v;\widetilde{K}_s) + 2\pi\theta_s(t)\big),
\end{align}
for every \(v\in\mathcal{A}_{\theta_s(t)}(\widetilde{K}_s)\); correspondingly, the argument appearing on the LHS runs over the interval \(\mathcal{A}_{\theta_t(s)}(\widetilde{K}_t)\).
\end{theorem}

\begin{proof}
Let \(E\) be any isoperimetric set of \(\widetilde K_s\) of volume
\(v\in\mathcal A_{\theta_s(t)}(\widetilde K_s)\). By
\cref{cor:iso-one-to-one} and \cref{rem:endpoints}, the set
\(
\frac{\delta_K(t)}{\delta_K(s)}\,E_{\theta_s(t)}
\)
is an isoperimetric set for \(\widetilde K_t\), of volume in
\(\mathcal A_{\theta_t(s)}(\widetilde K_t)\). Therefore,
\begin{equation}
\label{eq:TL-mid}
J\!\left(
\left|
\frac{\delta_K(t)}{\delta_K(s)}\,E_{\theta_s(t)}
\right|;
\widetilde K_t
\right)
=
\Per\!\left(
\frac{\delta_K(t)}{\delta_K(s)}\,E_{\theta_s(t)}
\right).
\end{equation}
Using the scaling properties of the area together with Steiner's area formula, which is valid since
\(|E|=v\in\mathcal A_{\theta_s(t)}(\widetilde K_s)\), see \cref{rem:use-At}, we obtain
\[
\left|
\frac{\delta_K(t)}{\delta_K(s)}\,E_{\theta_s(t)}
\right|
=
\frac{\delta_K(t)^2}{\delta_K(s)^2}
\Bigl(
|E|+\theta_s(t)\Per(E)+\pi\theta_s(t)^2
\Bigr)
=
\Phi_{s,t}(v).
\]
Similarly, using Steiner's perimeter formula,
\[
\Per\!\left(
\frac{\delta_K(t)}{\delta_K(s)}\,E_{\theta_s(t)}
\right)
=
\frac{\delta_K(t)}{\delta_K(s)}
\Bigl(
\Per(E)+2\pi\theta_s(t)
\Bigr).
\]
Since \(|E|=v\) and \(\Per(E)=J(v;\widetilde K_s)\), combining these identities with~\eqref{eq:TL-mid} yields~\eqref{eq:shifted-RL}.
\end{proof}

\begin{remark}
\label{rem:partialv}
Differentiating~\eqref{eq:shifted-RL} with respect to the spatial variable \(v\), for \(v\in\operatorname{int}\mathcal A_{\theta_s(t)}(\widetilde K_s)\), and setting
\(
w=\Phi_{s,t}(v),
\)
we obtain
\[
\partial_v\Phi_{s,t}(v)\,
\partial_wJ(w;\widetilde K_t)
=
\frac{\delta_K(t)}{\delta_K(s)}
\partial_vJ(v;\widetilde K_s).
\]
Using~\eqref{eq:Phi-derivative} and rearranging gives
\[
\partial_wJ(w;\widetilde K_t)
=
\frac{\delta_K(s)}{\delta_K(t)}
\cdot\frac{\partial_vJ(v;\widetilde K_s)}
{1+\theta_s(t)\,\partial_vJ(v;\widetilde K_s)}.
\]
This is consistent with~\eqref{eq:curvature-transformation} and with the fact that, for convex bodies, the derivative of the isoperimetric profile at a given volume equals the curvature of the free boundary of the corresponding isoperimetric set; see \cref{thm:LS22}~\ref{item:differentiability}. In particular, taking \(s=0\) recovers the corresponding relation between \(K\) and \(\widetilde K_t\).
\end{remark}


\subsection{Derivatives along the flow}
\label{sec:derivatives}
In this subsection, we compute the derivative of the reparametrization law~\eqref{eq:shifted-RL}.

\begin{theorem}[Derivative along the flow]
\label{thm:deriv-gen}
Let \(K\subset\R^2\) be a convex body. Then, for every \(w\in(0,|K|)\), the map
\(t\mapsto J(w;\widetilde K_t)\) is differentiable on \((-r_0,+\infty)\), and
\begin{equation}
\label{eq:transport_pde}
\frac{\mathrm{d}}{\mathrm{d}t} J(w;\widetilde K_t)\Big|_{t=s}
=
\delta_K(s)\,\mathcal D(w; \widetilde K_s),
\qquad \forall\, s>-r_0,
\end{equation}
where
\[
\mathcal D(w; L)
=
2\pi
-
\frac{\Per(L)}{2|L|}J(w;L)
-
\partial_w J(w;L)
\left(
J(w;L)-\frac{\Per(L)}{|L|}w
\right).
\]
\end{theorem}

\begin{proof}
Fix \(s>-r_0\) and \(w\in(0,|K|)\). By \cref{lem:Phi}~\ref{item:phi-IFT-domain}--\ref{item:phi-IFT-C1}, the map \(t\mapsto v(t)=\Phi_{s,t}^{-1}(w)\) is well defined and of class \(\mathrm C^1\) for \(|t-s|\ll1\), with \(v(s)=w\), and \(v(t)\) lies in the interior of \(\mathcal A_{\theta_s(t)}(\widetilde K_s)\). Hence, \cref{thm:RL} applies at \(v(t)\) and gives
\begin{equation}
\label{eq:shifted-at-w}
J(w;\widetilde K_t)
=
\frac{\delta_K(t)}{\delta_K(s)}
\Bigl(J(v(t);\widetilde K_s)+2\pi\theta_s(t)\Bigr),
\qquad |t-s|\ll1.
\end{equation}
By \cref{thm:LS22}~\ref{item:differentiability} and \cref{lem:area-C1}, the RHS is differentiable in \(t\), and so is the LHS.

We first compute \(v'(s)\).
Differentiating the identity \(\Phi_{s,t}(v(t))=w\) at \(t=s\), and recalling that \(\partial_v\Phi_{s,s}=1\) by \cref{lem:Phi}~\ref{item:phi-identity}, we get \(v'(s)=-\partial_t\Phi_{s,t}(w)\big|_{t=s}\). Since \(\theta_s(s)=0\) and \(\partial_t\theta_s(t)=\delta_K(s)\), differentiating~\eqref{eq:Phi} yields
\[
\partial_t\Phi_{s,t}(w)\Big|_{t=s}
=
2\frac{\delta_K'(s)}{\delta_K(s)}\,w+\delta_K(s)\,J(w;\widetilde K_s)
=
\delta_K(s)\left(J(w;\widetilde K_s)-\frac{\Per(\widetilde K_s)}{|K|}\,w\right),
\]
where we used~\eqref{eq:delta-derivative} in the form \(\delta_K'(s)/\delta_K(s)=-\delta_K(s)\Per(\widetilde K_s)/(2|K|)\).

We now differentiate~\eqref{eq:shifted-at-w} at \(t=s\). Using again~\eqref{eq:delta-derivative}, \(\theta_s(s)=0\), \(\partial_t\theta_s(t)=\delta_K(s)\) and \(v(s)=w\), we obtain
\[
\frac{\mathrm{d}}{\mathrm{d}t} J(w;\widetilde K_t)\Big|_{t=s}
=
-\frac{\delta_K(s)\Per(\widetilde K_s)}{2|K|}\,J(w;\widetilde K_s)
+\partial_wJ(w;\widetilde K_s)\,v'(s)
+2\pi\delta_K(s).
\]
Inserting the value of \(v'(s)\) computed above, and recalling that \(|\widetilde K_s|=|K|\), we get~\eqref{eq:transport_pde}.
\end{proof}

In the next proposition, we derive explicit characterizations of the derivative~\eqref{eq:transport_pde}, according to the range of the volume \(w\). We state it for \(t=0\) since it is the case we will use later on; nevertheless, the general case is stated in \cref{rem:derivative-RL-lv}.

\begin{proposition}[Explicit form of the derivative at \(t=0\)]
\label{prop:derivative-RL-lv}
Let \(K\subset\mathbb{R}^2\) be a convex body. Then: 

\noindent (i) for every \(w\in(0,\pi r_0^2]\), we have
\[
\frac{\mathrm{d}}{\mathrm{d}t}J(w;\widetilde K_t)\Big|_{t=0}=0;
\]

\noindent (ii) for every \(w\in(\pi r_0^2, \pi r_0^2 + 2r_0\mathcal{H}^1(K_{-r_0})]\), we have
\begin{equation}
\label{eq:TL-first-variation-large-vs-I}
\frac{\mathrm{d}}{\mathrm{d}t}J(w;\widetilde K_t)\Big|_{t=0}
= 
\ell \left(\frac{\Per(K)}{|K|} - \frac{2}{r_0} \right),
\end{equation}
where \(\ell \in (0, \mathcal{H}^1(K_{-r_0})]\) is the unique number such that \(w = \pi r_0^2 + 2r_0 \ell\);
\\
\noindent (iii) for every
\(w \in (\pi r_0^2 + 2r_0\mathcal{H}^1(K_{-r_0}), |K|)\),
we have
\begin{equation}
\label{eq:TL-first-variation-large-vs-II}
\hspace{-.2cm}
\frac{\mathrm{d}}{\mathrm{d}t}J(w;\widetilde K_t)\Big|_{t=0}
= 
\frac{\Per(K)}{r\delta_K(-r)^2} \left[1- \delta_K(-r) \left(1-\frac{\Per(K)}{2|K|}r \right)\frac{\Per(\widetilde{K}_{-r})}{\Per(K)}\right],
\end{equation}
where \(r\in[\rho_K,r_0)\) is the unique number such that \(w=|K_{-r}\oplus rB|\).
\end{proposition}

\begin{proof}
For case~(i), by \cref{thm:LS22}\ref{item:trivial-branch}, the isoperimetric set is a ball, and therefore \(J(w;K)=2\sqrt{\pi w}\). Inserting this expression into~\eqref{eq:transport_pde} immediately gives that the derivative vanishes. This is consistent with the enlargement of the ball branch of the profile following from Matheron's inequality~\eqref{eq:matheron}.

The remaining two cases are proved in the same way.

By \cref{thm:LS22}\ref{item:differentiability}, the volume-derivative \(\partial_v J(w;K)\) is the curvature of the free boundary of the isoperimetric set of volume \(w\). In case~(ii) this curvature equals \(1/r_0\), whereas in case~(iii) it equals \(1/r\). Depending on the range of \(w\), \cref{thm:LS22}\ref{item:first-nontrivial-branch} or \cref{thm:LS22}\ref{item:second-nontrivial-branch} describes the shape of the isoperimetric set: in case~(ii), it is the convex envelope of two maximal balls whose centers are separated by a distance \(\ell\); in case~(iii), it is \(K_{-r}\oplus rB\) for a unique \(r\in[\rho_K,r_0)\). In both cases, one can compute the perimeter and volume explicitly, and the claimed formulas follow by direct algebraic simplification.
\end{proof}

\begin{remark}
\label{rem:derivative-RL-lv}
By \cref{prop:derivative-RL-lv} applied to \(\widetilde{K}_s\) together with \cref{thm:deriv-gen}, we can write the derivative at time \(t=s>-r_0\) as follows.
For \(w\in(0,\pi r_s^2]\), the derivative vanishes. If \(w\in (\pi r_s^2,\pi r_s^2+2r_s(\mathcal{H}^1({(\widetilde{K}_s)}_{-r_s}))]\), then
\[
\frac{\mathrm{d}}{\mathrm{d}t} J(w;\widetilde K_t)\Big|_{t=s}
=
\delta_K(s)\ell \left(\frac{\Per(\widetilde{K}_s)}{|\widetilde{K}_s|}-\frac{2}{r_s}\right),
\]
where \(\ell\in (0, \mathcal{H}^1((\widetilde{K}_s)_{-r_s})]\) is the unique number such that \(w=\pi r_s^2 + 2r_s\ell\).

If instead, \(w>\pi r_s^2+2r_s\mathcal{H}^1({(\widetilde{K}_s)}_{-r_s})\), then
\begin{equation*}
\frac{\mathrm{d}}{\mathrm{d}t} J(w;\widetilde K_t)\Big|_{t=s}
=
\delta_K(s) \Delta(w; \widetilde{K}_s),
\end{equation*}
where
\[
\Delta(w; L)
=
\frac{\Per(L)}
{r\,\delta_{L}(-r)^2}
\left[
1
-
\delta_{L}(-r)
\left(
1-\frac{\Per(L)}{2|L|}r
\right)
\frac{\Per(
{\widetilde{L}}_{-r}
)}{\Per(L)}
\right],
\]
with
\(r=r(w)\in
[\rho_L, \inr(L))
\)
the unique number such that \(w=|L_{-r}\oplus rB|\).
\end{remark}


\subsection{Consequences of the reparametrization law}
\label{ssec:consequence-RL}


\subsubsection{Comparison on different bodies}
\label{ssec:comparison}

While \cref{thm:monotonicity} allows us to compare the profile of a body \(K\) with its evolved set \(\widetilde{K}_t\), as a consequence of~\eqref{eq:shifted-RL}, we are able to show that if two different bodies (of equal area), \(K, L\) are such that their profile satisfy
\[
J(v;K)\ge J(v;L), \qquad \forall v,
\]
then, the above inequality is preserved along the flow.

\begin{theorem}[Order preservation along the flow]
\label{thm:comparison}
Let \(K,L\subset\R^2\) be convex bodies of same area, and assume that \(J(v;K)\ge J(v;L)\) for every \(v\in[0,|K|]\). Then, the inequality is preserved along the  (forward-in-time) flow, that is,
\[
J(v;\widetilde K_t)\ge J(v;\widetilde L_t)
\]
for every \(v\in[0,|K|]\) and every \(t\ge0\).
\end{theorem}

\begin{proof}
Evaluating the assumption at \(v=|K|\) gives \(\Per(K)\ge\Per(L)\), so that, by Steiner's formula~\eqref{eq:steiner-area}, \(|K_t|\ge|L_t|\) and hence
\(\delta_K(t)\le\delta_L(t)\) for every \(t\ge0\).

\smallskip
\underline{\emph{Step~(i). Unnormalized profiles.}}
For the ease of notation, let \(f, f_t\) denote, respectively, \(J(-; K)\) and \(J(-;K_t)\), while \(g, g_t\) the profiles with respect the ambient bodies \(L\) and \(L_t\).

Fix \(t\ge0\). Since every volume is admissible, \cref{thm:RL} with \(s=0\), together with the rescaling of \cref{lem:scaling-profile}, gives
\begin{equation}
\label{eq:RL-unnormalized}
f_t(v+t f(v)+\pi t^2)=f(v)+2\pi t,
\qquad v\in[0,|K|],
\end{equation}
and the same for \(L\) with \(g,g_t\) in place of \(f, f_t\). By \cref{thm:LS22}~\ref{item:differentiability}, all functions \(f,f_t,g,g_t\) are continuous and strictly increasing. Letting \(y=f(v)+2\pi t\), using that \(f,f_t\) are invertible, from~\eqref{eq:RL-unnormalized}, we find
\begin{equation}
\label{eq:inverse}
f^{-1}_t (y)
=
f^{-1}(y-2\pi t) + t(y-2\pi t) + \pi t^2,
\qquad
y\in [2\pi t, \Per(K_t)],
\end{equation}
and the same holds for the inverse functions of \(g,g_t\) on \([2\pi t, \Per(L_t)]\).

Notice that the standing assumption \(f(v)\ge g(v)\) on \([0,|K|]\) is equivalent to \(f^{-1}(z) \le g^{-1}(z)\) on \([0, \Per(L)]\).

By~\eqref{eq:inverse} on \(f,f_t\), then using \(f^{-1}(z) \le g^{-1}(z)\), and once again by~\eqref{eq:inverse} on \(g,g_t\), we find
\[
f^{-1}_t(y) \le g^{-1}_t(y),
\qquad
\forall y \in [2\pi t, \Per(L_t)].
\]
This latter is once again equivalent to the opposite inequality for the functions \(f_t,g_t\), so that
\begin{equation}
\label{eq:ft-to-gt}
f_t(v) \ge g_t(v),
\qquad
\forall v\in [\pi t^2, |L_t|].
\end{equation}
For \(v\le\pi t^2\), both profiles equal \(2\sqrt{\pi v}\), since the inradii of \(K_t\) and \(L_t\) are at least \(t\). Thus \(f_t(v)\ge g_t(v)\) for every \(v\in[0,|L_t|]\).

\smallskip
\underline{\emph{Step~(ii). Conclusion.}}
We now transpose the inequality of the previous step to the rescaled profiles. Fix \(v\in[0,|K|]\), and, up to a translation, suppose that \(0\in K_t\). Since \(\delta_K(t)\le\delta_L(t)\), convexity gives \(\widetilde K_t=\delta_K(t)K_t\subseteq\delta_L(t)K_t\), so that for all \(v\le |\widetilde{K}_t|\) the profile \(J(-; \widetilde{K}_t)\) is bound from below by that of \(J(-; \delta_L(t)K_t)\).

By the rescaling of the profile of \cref{lem:scaling-profile}, Step~(i), and the rescaling once again, we get
\begin{align*}
J(v; \widetilde{K}_t) 
\ge
J(v; \delta_L(t)K_t)
&
=
\delta_L(t) J\left(\frac{v}{\delta_L(t)^2}; K_t\right)
\\
&
\ge
\delta_L(t) J\left(\frac{v}{\delta_L(t)^2}; L_t\right)
=
J(v; \widetilde{L}_t).
\qedhere
\end{align*}
\end{proof}


\subsubsection{\texorpdfstring{\(\mathrm{C}^1\)}{C1}-regularity of the flow and Hamilton--Jacobi equation}
\label{ssec:HJ}
A first consequence of the reparametrization law proved in \cref{thm:RL}, together with \cref{lem:Phi}, is the \(\mathrm{C}^1\) continuity of the isoperimetric profile along the flow in the joint couple \((v,t)\).

\begin{theorem}[Regularity of the flow]
\label{thm:reg-flow}
Let \(K\subset\mathbb R^2\) be a convex body. Then, the map
\[
(w,t)\mapsto J(w;\widetilde K_t)
\]
is continuous on \([0,|K|]\times(-r_0,+\infty)\) and of class \(\mathrm{C}^1\) on the interior.
\end{theorem}
\begin{proof}
Let \((w_0, t_0)\) be fixed in \((0,|K|)\times(-r_0,+\infty)\) and set \(s=t_0\).
By \cref{lem:Phi}~\ref{item:phi-IFT-domain} we can use the shifted reparametrization law. Hence,
\[
J(w;\widetilde K_t)
=
\frac{\delta_K(t)}{\delta_K(t_0)}
\Bigl(
J(\Phi^{-1}_{t_0,t}(w);\widetilde K_{t_0})+2\pi\theta_{t_0}(t)
\Bigr),
\]
where \((w,t)\mapsto \Phi^{-1}_{t_0,t}(w)\) is of class \(\mathrm{C}^1\) by \cref{lem:Phi}~\ref{item:phi-IFT-C1}. By \cref{thm:LS22}~\ref{item:differentiability}, the definition of \(\theta_{t_0}(t)\), and \cref{lem:area-C1}, the RHS is of class \(\mathrm C^1\).

To conclude the proof, we have to show the continuity of the flow when \(w_0\) is either \(0\) or \(|K|\).

\underline{\emph{Case~(i): \(w_0=0\).}}
Clearly, \(w_0<\pi r_{t_0}^2\).
By continuity of \(t\mapsto r_t\), there exists a neighborhood of \((w_0,t_0)\) in which \(w<\pi r_t^2\). Hence, by \cref{thm:LS22}~\ref{item:trivial-branch},
\(
J(w;\widetilde K_t)=2\sqrt{\pi w}
\)
throughout this neighborhood, and continuity follows.

\underline{\emph{Case~(ii): \(w_0=|K|\).}}
Fix \(\varepsilon>0\) and choose
\(w_*\in(0,|K|)\) such that
\begin{equation}
\label{eq:choice-wstar}
J(|K|;\widetilde K_{t_0})-J(w_*;\widetilde K_{t_0})<\varepsilon,
\end{equation}
which is possible by continuity of \(J(-;\widetilde K_{t_0})\), see \cref{thm:LS22}~\ref{item:differentiability}. 

Since \(w_*\in(0,|K|)\), the map \(t\mapsto J(w_*;\widetilde K_t)\) is continuous at \(t_0\) by the first part of the proof; moreover, \(t\mapsto J(|K|;\widetilde K_t)=\delta_K(t)\Per(K_t)\) is continuous by \cref{lem:area-C1}. Hence, for \(|t-t_0|\ll 1\),
\begin{equation}
\label{eq:two-continuities}
|J(w_*;\widetilde K_t)-J(w_*;\widetilde K_{t_0})|<\varepsilon
\quad\text{and}\quad
|J(|K|;\widetilde K_t)-J(|K|;\widetilde K_{t_0})|<\varepsilon .
\end{equation}

Let now \(w\in[w_*,|K|]\) be fixed and \(|t-t_0|\ll 1\). 
On the one hand, we can write the difference \(J(|K|;\widetilde K_{t_0})-J(w;\widetilde K_t)\) as
\[
\bigl[J(|K|;\widetilde K_{t_0})-J(w_*;\widetilde K_{t_0})\bigr]
+
\bigl[J(w_*;\widetilde K_{t_0})-J(w_*;\widetilde K_t)\bigr]
+
\bigl[J(w_*;\widetilde K_t)-J(w;\widetilde K_t)\bigr].
\]
The first bracket is smaller than \(\varepsilon\) by~\eqref{eq:choice-wstar}, the second is smaller than \(\varepsilon\) by~\eqref{eq:two-continuities}, and the third is non-positive since \(w\ge w_*\) and \(J(-;\widetilde{K}_t)\) is increasing. Hence,
\(J(|K|;\widetilde K_{t_0})-J(w;\widetilde K_t)<2\varepsilon\).

On the other hand, since \(w\le|K|\), owing again to the fact that \(J(-;\widetilde{K}_t)\) is increasing, we have
\[
J(w;\widetilde K_t)-J(|K|;\widetilde K_{t_0})
\le
J(|K|;\widetilde K_t)-J(|K|;\widetilde K_{t_0})
\le
\varepsilon,
\]
by~\eqref{eq:two-continuities}. This proves the claimed continuity.
\end{proof}

As a consequence of \cref{thm:deriv-gen} and of \cref{thm:reg-flow}, the following holds.

\begin{theorem}
\label{thm:HJ_PDE}
Let \(K\subset\mathbb{R}^2\) be a convex body, \(p(t) = \Per(\widetilde{K}_t)\), and \(\delta_K(t)\) be as in~\eqref{eq:renormalization}. 
The isoperimetric profile along the volume-preserving Minkowski regularizing flow \(t\mapsto\widetilde{K}_t\) is a classical solution on \((0,|K|)\times(-r_0,+\infty)\) of the first-order quasilinear Hamilton--Jacobi equation
\begin{equation}
\label{eq:HJ-PDE}
\tag{TL}
\partial_t I(w,t)
=
\delta_K(t)
\left(
2\pi - \frac{p(t)}{2|K|}I(w,t) - \partial_w I(w,t)\left(I(w,t)-\frac{p(t)}{|K|}w \right)
\right),
\end{equation}
it is continuous up to \(w\in\{0,|K|\}\), and it satisfies the lateral
boundary conditions
\[
I(0,t)=0,
\qquad
I(|K|,t)=p(t),
\qquad
\forall\, t>-r_0 .
\]
In particular, its restriction to \((0,|K|)\times(0,+\infty)\) is a classical solution of the Cauchy problem for~\eqref{eq:HJ-PDE} with initial datum \(I(w,0)=J(w;K)\).
\end{theorem}

\begin{remark}[Dimensional analysis]
Although $t$ is conventionally referred to as ``time'', it represents the radius of the regularizing ball and carries units of length \([\unit{m}]\). The state variable $w$ has units of \([\unit{m^2}]\), while the isoperimetric profile \(J(w;\widetilde{K}_t)=I(w,t)\) has units of length \([\unit{m}]\). Consequently, the transport law~\eqref{eq:HJ-PDE} is dimensionally consistent.
\end{remark}

The PDE~\eqref{eq:HJ-PDE} will be the object of a forthcoming paper. Here, we just make the following observations.

\begin{remark}[Lateral conditions]
Since~\eqref{eq:HJ-PDE} is of first order, it transports information along curves \(w=w(t)\) whose velocity is the coefficient of \(\partial_wI\). 
Evaluated along the solution, this velocity vanishes at both \(w=0\) and \(w=|K|\), so that the lateral boundary is characteristic. 
The lateral conditions should therefore not be regarded as data to be prescribed, but as compatibility conditions satisfied by the solution.
\end{remark}

\begin{remark}[Uniqueness and comparison]
\label{rem:HJ-perspectives}
Equation~\eqref{eq:HJ-PDE} is not a standard Hamilton--Jacobi equation, and the usual comparison principles for viscosity solutions do not apply directly. First, it is nonlocal: \(p(t)=I(|K|,t)\) is the trace of the solution itself, and \(\delta_K\) is in turn determined by \(p\) through \cref{lem:area-C1}. 
Second, the coupling \(\partial_wI\cdot I\) makes the Hamiltonian proper in \(I\) along non-decreasing profiles, but not in general.
Third, the lateral boundary is
characteristic and singular, see \cref{rem:boundary-degeneracy}.
We stress that \cref{thm:monotonicity} is not a comparison statement, but a sign condition on \(\partial_tI\), obtained by variational means. Whether~\eqref{eq:HJ-PDE} characterizes the flow of profiles remains open. Nonetheless, the ordering of the profiles of bodies with equal area is preserved along the flow for \(t\ge0\), see \cref{thm:comparison}; the proof relies on the reparametrization law rather than on a maximum principle.
\end{remark}

\begin{remark}[Unnormalized profiles]
\label{rem:burgers}
For \(t\ge0\), set \(\hat I(w,t)=J(w;K_t)\), \(w\in[0,|K_t|]\), which is of class \(\mathrm C^1\) on \((0,|K_t|)\times(0,+\infty)\) by the rescaling of the profiles of  \cref{lem:scaling-profile} and the regularity of non-rescaled profiles of \cref{thm:reg-flow}. By the reparametrization law~\eqref{eq:shifted-RL} and the scaling \cref{lem:scaling-profile} (see also~\eqref{eq:RL-unnormalized}), along each curve \(w_v(t)=v+tJ(v;K)+\pi t^2\), with \(v\in[0,|K|]\), one has
\[
\dot w_v(t)
=
J(v;K)+2\pi t
=
\hat I\bigl(w_v(t),t\bigr),
\qquad
\frac{\mathrm d}{\mathrm dt}\hat I\bigl(w_v(t),t\bigr)
=
2\pi.
\]
Thus, the volume moves with speed equal to the value of the profile at the moving point, while the latter increases at rate \(2\pi\). For \(t\ge0\), the map \(v\mapsto w_v(t)\) is strictly increasing, with derivative \(1+t\,r_K(v)^{-1}>0\), and maps \([0,|K|]\) onto \([\pi t^2,|K_t|]\); hence,
by the chain rule,
\[
\partial_t\hat I+\hat I\,\partial_w\hat I
=
2\pi
\]
on \(\{\pi t^2<w<|K_t|\}\). The same equation holds on \(\{0<w<\pi t^2\}\), where \(\hat I(w,t)=2\sqrt{\pi w}\); these points lie on the curves \(w=\pi(t-t_0)^2\), \(t_0>0\), which leave \(w=0\) at positive times, in accordance with the singularity of \(\partial_w\hat I\) at \(w=0\). Thus \(\hat I\) is a classical solution of Burgers' equation with constant forcing, which is local and the same for every body: the nonlocality of~\eqref{eq:HJ-PDE} is entirely due to the renormalization~\eqref{eq:renormalization}. The right boundary \(w=|K_t|\) is the curve \(w_{|K|}\), since by \cref{lem:area-C1} \(\frac{\mathrm{d}}{\mathrm{d}t}|K_t|=\Per(K_t) =\hat I(|K_t|,t)\). Finally, for \(t<0\) the derivative of the map \(v\mapsto w_v(t)\) is strictly positive exactly where \(r_K(v)>-t\), which is the admissibility condition~\eqref{eq:admissible} up to its endpoints: for negative times, the loss of admissibility corresponds to the crossing of the curves \(w_v\).
\end{remark}


\subsubsection{Concavity of the volume of inner parallel sets}
\label{sec:concavity}
In this subsection, we turn our attention to \cref{prop:derivative-RL-lv}.

If \(K\) has a unique inball, case~(ii) of \cref{prop:derivative-RL-lv} never occurs. Assuming instead that \(K_{-r_0}\) is not a singleton, the derivative in~\eqref{eq:TL-first-variation-large-vs-I} is always negative by the classical inequality relating area, perimeter, and inradius for convex sets, see~\cite[Tab.~3.1, triplet \((A,p,r)\)]{SA00}.

We now consider equation~\eqref{eq:TL-first-variation-large-vs-II}. By \cref{thm:monotonicity}, this quantity must be non-positive. This is equivalent to
\begin{equation}
\label{eq:derivative-upper-branch-sign}
\delta_K(-r)\left(1-\frac{\Per(K)}{2|K|}r\right)
\left(\frac{\Per(\widetilde K_{-r})}{\Per(K)}\right)\ge 1.
\end{equation}
By \cref{prop:rs}, the last factor is at least \(1\), with equality only for tangential bodies. The product of the first two factors is also at least \(1\) as a direct consequence of the concavity of \(x\mapsto |K_x|^{\frac12}\) stemming from properties of parallel bodies~\cite[Lem.~3.1.13]{Sch14book} and the classical Brunn--Minkowski inequality~\cite[Thm.~7.1.1]{Sch14book}. Exploiting \cref{prop:rs}, see also~\cite[Thm.~1.1]{RS22}, we can characterize cases of equality.

\begin{proposition}
\label{prop:ilias-ineq}
Let \(K\subset\mathbb{R}^2\) be a convex body. Then, for all 
\(t,s>-r_0\),
\begin{equation}
\label{eq:ilias-ineq}
|K_{t}|^{\frac 12}
\le
|K_s|^{\frac12}\left(1+\frac{\Per(K_s)}{2|K_s|}(t-s)\right).
\end{equation}
In particular, equality is attained at \(t\neq s\) if and only if \(\max\{t,s\}\le \tau_K\).
\end{proposition}

\begin{proof}
By~\cite[Lem.~3.1.13 and Thm.~7.1.1]{Sch14book} the function \(\psi:(-r_0,\infty)\to\mathbb{R}\) defined by
\(
\psi(x)=|K_x|^{\frac12}
\)
is concave, therefore it rests below its tangent lines. Recalling \cref{lem:area-C1}, its derivative is
\begin{equation}
\label{eq:psi-derivative}
\psi'(x)
=
\frac12\,\frac{\Per(K_x)}{|K_x|^{\frac12}}
=
\frac{1}{{2|K|^{\frac12}}}\,\Per(\widetilde K_x),
\end{equation}
so that~\eqref{eq:ilias-ineq} follows.

For equality, note that if \(t\neq s\) and equality holds, then \(\psi\) must be affine on the interval joining \(s\) and \(t\). By~\eqref{eq:psi-derivative}, this happens if and only if \(\Per(\widetilde K_x)\) is constant in that interval, which by \cref{prop:rs} occurs precisely when \(x\in(-r_0,\tau_K]\). Hence, equality for some \(t\neq s\) holds if and only if \(t,s\le\tau_K\).
\end{proof}

\begin{remark}
\label{rem:ilias-ineq}
Exploiting the \(N\)-dimensional analogue of \cref{prop:rs} proved in~\cite[Thm.~1.1]{RS22}, one obtains the \(N\)-dimensional version
\[
|K_t|^{\frac1N}
\le
|K_s|^{\frac1N}\left(1+\frac{\Per(K_s)}{N|K_s|}(t-s)\right)
\]
by the same argument.
\end{remark}

As a consequence of \cref{prop:ilias-ineq}, we are able to provide a new proof of~\cite[Eq.~(5)]{Fto23}. We also completely characterize cases of equality, completing the partial picture provided in \cite{Fto23}.

\begin{corollary}
\label{cor:ilias}
Let \(K\subset\mathbb{R}^2\) be a convex body and let \(h_1(K)\) be its Cheeger constant. It holds
\[
h_1(K)
\ge 
\frac{\Per(K)+\sqrt{4\pi |K|}}{2|K|}.
\]
Equality holds if and only if $K$ is a tangential body. 
\end{corollary}

\begin{proof}
It is well known that if \(K\) is a tangential body, then equality holds, see~\cite[Eq.~(16)]{Fto23}. Let us now assume that \(K\) is not a tangential body, that is, \(\tau_K<0\). By~\eqref{eq:ilias-ineq}, for every \(t\in (0,r_0)\), we have
\begin{equation}
\label{eq:ineq_str}
|K|\left(1-\frac{\Per(K)}{2|K|}t\right)^2
>
|K_{-t}|. 
\end{equation}
Let us denote by \(g(t)\) the LHS, as a function of \(t\).

It is well known by~\cite{KLR06} (see also~\cite[Cor.~5.5]{LS20} for an extension to non-convex sets) that the Cheeger constant of $K$ is the inverse of \(t^*\), the (unique) solution of equation \(|K_{-t}|=\pi t^2\) in \((0,r_0)\). By~\eqref{eq:ineq_str}, the zero \(t^*=1/h_1(K)\) is  strictly bounded from above by the smallest positive solution of the equation \(g(t)=\pi t^2\), given by
\[
\frac{2|K|}{\Per(K)+\sqrt{4\pi |K|}}.
\]
Therefore, the strict inequality follows when \(K\) is not a tangential body.
\end{proof}


\subsubsection{Profile differential inequality}

As a consequence of \cref{thm:monotonicity} and of \cref{thm:deriv-gen}, we derive the following differential inequality.

\begin{theorem}
\label{thm:differential_ineq}
Let \(K\subset\mathbb{R}^2\) be a convex body. Then, for all \(v\in (0,|K|)\), the isoperimetric profile \(J(v;K)\) satisfies the differential inequality
\begin{equation}
\label{eq:differential_ineq}
2\pi - \frac{\Per(K)}{2|K|}J(v;K) - \partial_v J(v; K)\left(J(v;K)-\frac{\Per(K)}{|K|}v \right)
\le
0.
\end{equation}
In particular, this is an equality for all \(v\le \pi \inr(K)^2\), whereas it is an equality for some (equivalently, every) \(v\in(\pi \inr(K)^2,|K|)\) if and only if \(K\) is a tangential body.
\end{theorem}

\begin{proof}
The inequality follows directly from \cref{thm:monotonicity} and \cref{thm:deriv-gen}.

Let us prove equality cases. If \(v\le \pi\inr(K)^2\), equality follows from \cref{prop:derivative-RL-lv}~(i). Assume now \(v>\pi\inr(K)^2\). If \(v\) lies in the range of \cref{prop:derivative-RL-lv}~(ii), this latter is non-empty. Then \(K\) has more than one inball, and hence it is not a tangential body. In particular, the derivative reads as in~\eqref{eq:TL-first-variation-large-vs-I}, which is strictly negative for non-tangential bodies. Otherwise, \(v\) lies in the range of \cref{prop:derivative-RL-lv}~(iii), and the LHS of~\eqref{eq:differential_ineq} vanishes if and only if equality holds in~\eqref{eq:derivative-upper-branch-sign}, which by \cref{prop:rs} and \cref{prop:ilias-ineq} happens if and only if \(K\) is a tangential body.
\end{proof}
From the above differential inequality, we can prove the following sharp lower bound to the isoperimetric profile \(J(v;K)\), uniquely attained by tangential bodies.
\begin{proposition}
\label{prop:lower-bound-iso-profile}
Let \(K\subset\mathbb R^2\) be a convex body and let \(v\in[0,|K|]\). Then,
\[
J(v;K)
\ge
\begin{cases}
2\sqrt{\pi v},
&\text{if } v\le v^*,\\[6pt]
\Per(K)
-\sqrt{\dfrac{\Per(K)^2-4\pi|K|}{|K|}\,(|K|-v)},
&\text{otherwise},
\end{cases}
\]
where
\[
v^*
=
\frac{4\pi |K|^2}{\Per(K)^2}.
\]
Equality on the whole range \([0,|K|]\) is uniquely attained by tangential bodies.
\end{proposition}

Clearly, the first bound comes from the isoperimetric inequality. By~\cite[Table~3.1, triplet \((A,p,r)\)]{SA00}, \(v^*\) is at least \(\pi \inr(K)^2\). Indeed,
\[
v^*
=
\frac{4\pi|K|^2}{\Per(K)^2}
=
\frac{\pi (2|K|)^2}{\Per(K)^2} 
\ge
\pi \inr(K)^2.
\]
Then, in the second regime of the above proposition, that is, when \(v>v^*\), we are in the convex regime of \cref{thm:LS22} (see points~\ref{item:first-nontrivial-branch}--\ref{item:second-nontrivial-branch} and~\ref{item:convexity}).  This lower bound is strictly convex in \(v\), depends only on \(|K|\) and \(\Per(K)\), and is sharp. By contrast, the linear lower bounds provided by the convexity of \(J(-;K)\), namely its tangent lines, require the knowledge of \(J\) and of its derivative at some volume, such as the inradius or the rolling radius of \(K\).

\begin{proof}
For the sake of brevity, let \(A=|K|\) and \(P=\Per(K)\).

By the isoperimetric inequality, the lower bound \(2\sqrt{\pi v}\) always holds. If \(K\) is a ball, then \(P^2=4\pi A\) and there is nothing to prove. 
Thus, for the remainder of the proof, we assume
\[
P^2>4\pi A
\qquad
\text{and }
v > \frac{4\pi A^2}{P^2}.
\]
By the isoperimetric inequality, we have
\begin{equation}
\label{eq:rough-est}
J(v;K) 
\ge 
2\sqrt{\pi v}
>
\frac{4\pi A}{P}.
\end{equation}

Let us consider the auxiliary function defined for \(v\in (4\pi A^2/P^2, A]\)
\[
\Phi(v)
=
\frac{
v-\dfrac{2A}{P}J(v;K)
+\dfrac{4\pi A^2}{P^2}
}
{
\left(
2\pi-\dfrac{P}{2A}J(v;K)
\right)^2
},
\]
whose derivative is
\[
\Phi'(v)
=
\frac{
2\pi-\dfrac{P}{2A}J(v;K)
-\partial_vJ(v;K)
\left(
J(v;K)-\dfrac{P}{A}v
\right)
}{
\left(
2\pi-\dfrac{P}{2A}J(v;K)
\right)^3
}.
\]
By \cref{thm:differential_ineq} and~\eqref{eq:rough-est}, we have \(\Phi'\ge 0\), thus,
\begin{equation}
\label{eq:final1}
\Phi(v) 
\le 
\Phi(A) 
=
-
\frac{4A^3}{P^2(P^2-4\pi A)}.
\end{equation}
Writing explicitly \(\Phi(v)\), multiplying by the positive denominators, simplifying, and rearranging, we get
\begin{equation}
\label{eq:final2}
A(P-J(v;K))^2 \le (P^2-4\pi A)(A-v),
\end{equation}
from which the claimed lower bound follows.

Let us turn our attention to the equality cases. We notice that, by~\cite[Table~3.1, triplet \((A,p,r)\)]{SA00}, \(v\) is at least \(\pi \inr(K)^2\).
Hence, by \cref{thm:differential_ineq}, \(\Phi'(v)=0\) if and only if \(K\) is a tangential body. Therefore, inequalities~\eqref{eq:final1}--\eqref{eq:final2} become equalities. Since for tangential bodies it holds \(\inr(K)=2A/P\), the lower bound \(2\sqrt{\pi v}\) is saturated for volumes smaller than \(4\pi A^2/P^2\).
\end{proof}


\subsection{Two degeneracies}
\label{sec:degeneracies}
In this subsection, we briefly describe the two phenomena that limit the local, differential perspective. The first one explains why the global proof of \cref{thm:monotonicity} cannot be dispensed with; the second one describes what happens below the rolling radius \(\rho_K\), where the flow is rougher, in a way which is invisible in the formulas of the previous subsections.


\subsubsection{Degeneracy at tangential bodies}
\label{ssec:degeneracy-tangential}
The local behavior of \(J(w;\widetilde K_t)\) at \(t=0\), together with Matheron's inequality~\eqref{eq:matheron}, Richter--Saor\'in-G\'omez proof of the inequality in \cref{prop:rs}, and the inequality in \cref{prop:ilias-ineq} are not sufficient to prove \cref{thm:monotonicity}. Indeed, as stated in \cref{thm:differential_ineq}, the derivative formula degenerates for tangential bodies, at every volume. Tangential bodies are exactly the stationary points of the flow, by~\eqref{eq:constancy-interval}, and on them one would need \emph{at least} a second-order infinitesimal information to deduce monotonicity. The global argument completely avoids the need to compute higher-order derivatives.


\subsubsection{Degeneracy below the rolling radius}
\label{ssec:degeneracy-rho}
\cref{lem:outer-composition} gives the composition law for non-negative \(s,t\). If the rolling radius of \(K\) is positive, \(\rho_K>0\), \cref{lem:outer-composition} extends verbatim to \(t>s>-\rho_K\), since it still holds \(K_t = (K_s)_{t-s}\). On the other hand, if \(s<-\rho_K\) it fails, since, e.g., \((K_s)_{-s}\) is a proper subset of \(K\).

In particular, when restricting to this range, all proofs simplify as one can exploit the composition law for prefactors contained in \cref{lem:outer-composition}, so that from statements on the pair given by the body \(K\) and any evolved set \(\widetilde{K}_t\), one can derive that for general pairs.

For \(t\le -\rho_K\), the evolved set \(\widetilde{K}_t\) loses \(\mathrm{C}^{1,1}\) boundary regularity.  This has some consequences at the level of the Hamilton--Jacobi
equation~\eqref{eq:HJ-PDE}. First, the solution enjoys extra regularity precisely for \(t>-\rho_K\), as the next remark shows.

\begin{remark}[Degeneracy of the spatial derivative]
\label{rem:boundary-degeneracy}
By \cref{thm:HJ_PDE}, the flow classically solves the quasilinear Hamilton--Jacobi PDE~\eqref{eq:HJ-PDE}.
By~\eqref{eq:RLv} and \cref{thm:LS22}~\ref{item:differentiability-lv}, one has \(\partial_wI(w,t)=r_{\widetilde K_t}(w)^{-1}\), and \(r_{\widetilde K_t}(w)\to\rho_{\widetilde K_t}\) as \(w\to|K|\). Therefore,
\[
\lim_{w\to|K|}\partial_wI(w,t)<+\infty
\qquad\Longleftrightarrow\qquad
\rho_{\widetilde K_t}>0.
\]
The latter holds true if and only if \(t>-\rho_K\).
Indeed, since \(\widetilde K_t=\delta_K(t)K_t\), one has \(\rho_{\widetilde K_t}=\delta_K(t)\rho_{K_t}\), and hence \(\rho_{\widetilde K_t}>0\) if and only if \(\rho_{K_t}>0\).

If \(t>-\rho_K\), then \(\sigma_{K_t}(u)=\sigma_K(u)+t\) for every \(u\in\mathbb S^1\). Hence, by the definition of the rolling radius, \(\rho_{K_t}=\rho_K+t>0\).
Conversely, assume that \(\rho_{K_t}>0\). The case \(t\ge0\) is immediate, so let \(t<0\). Since \(\rho_{K_t}>0\), \(K_t\) is regular, that is, it has a unique supporting line at each boundary point. The argument in the proof of~\cite[Lem.~2.5]{HS14} yields \(\sigma_{K_t}(u)=\sigma_K(u)+t\) for every \(u\in\mathbb S^1\). Therefore \(\rho_K=\rho_{K_t}+|t|>|t|\), and hence \(t>-\rho_K\).

Thus, the flow extends in a \(\mathrm{C}^1\) fashion at \(w=|K|\) if and only if \(t>-\rho_K\).
\end{remark}

Second, the boundary datum \(p(t)=\Per(\widetilde{K}_t)\), which also appears as a coefficient, is not, in general, of class \(\mathrm{C}^1\) in any open neighborhood of \(t=-\rho_K\).

\begin{remark}
\label{rem:lipschitz-p}
Steiner's perimeter formula~\eqref{eq:steiner-per} holds for \(t>-\rho_K\); therefore, the boundary datum \(p(t)=\Per(\widetilde{K}_t)\) is of class \(\mathrm{C}^1\) in this range.

On the whole range, \(p\) is locally Lipschitz. Indeed, let \(-r_0<s<t\le0\). Then \(K_s\) is the inner parallel body of \(K_t\) at distance \(t-s\), and \(\inr(K_t)=r_0+t\). Larson's inequality~\cite{Lar16} applied to \(K_t\), together with the monotonicity of the perimeter under inclusion of convex sets, yields
\[
0\le \Per(K_t)-\Per(K_s)\le\frac{\Per(K_t)}{r_0+t}\,(t-s).
\]
Therefore, \(t\mapsto\Per(K_t)\) is locally Lipschitz on \((-r_0,0]\), and so is \(p\), since \(\delta_K\) is of class \(\mathrm C^1\) by \cref{lem:area-C1}.
In general, one cannot expect \(\mathrm{C}^1\) regularity in any neighborhood of \(t=0\), as polygons show. 

Let \(K\) be a convex polygon, so that \(\rho_K=0\). Let us denote its exterior angles by \(\alpha_1,\dots,\alpha_m\in(0,\pi)\), which sum up to \(2\pi\). For \(t\ge0\), by Steiner's perimeter formula~\eqref{eq:steiner-per}, we have \(\Per(K_t)=\Per(K)+2\pi t\). For \(t\le0\) with \(|t|\) small, each edge of \(K_t\) is shortened by \(|t|\tan(\alpha_j/2)\) at the vertex of exterior angle \(\alpha_j\), so that \(\Per(K_t)=\Per(K)+2t\sum_j\tan(\alpha_j/2)\). Since \(\delta_K\) is of class \(\mathrm C^1\), the one-sided derivatives of \(p\) at
\(t=0\) satisfy
\[
p'(0^-)-p'(0^+)
=
2\sum_j\tan(\alpha_j/2)-2\pi>0,
\]
because \(\tan x>x\) on \((0,\pi/2)\).
For instance, if \(K\) is a tangential polygon, then \(p\) is constant on \((-r_0,0]\) by~\eqref{eq:constancy-interval}, whereas \(p'(0^+)=2\pi-\Per(K)/r_0<0\), since \(\Per(K)=2r_0\sum_j\tan(\alpha_j/2)>2\pi r_0\). This is consistent with \cref{thm:differential_ineq}, which gives \(\partial_tI(w,0)=0\) for every \(w\in(0,|K|)\): at \(w=|K|\), where \(\partial_wI\) blows up, see \cref{rem:boundary-degeneracy}, the time derivative and the limit \(w\to|K|\) do not commute.
\end{remark}


\section{Non-convex sets}
\label{sec:non-convex}

It is natural to ask whether analogous results hold for non-convex sets. The main tools used above come from~\cite{LS22}, which applies to a broader class of (open) planar sets introduced in~\cite{LNS17}. We recall that class here.

\begin{definition}[Set with no-necks]
\label{def:no-neck}
We say that a set \(\Omega\subset\mathbb{R}^2\) has \emph{no-necks} if the following conditions hold:
\begin{enumerate}[label=(\roman*)]
\item
\label{item:no-necks-regularity}
\(\Omega\) is a Jordan domain and the two-dimensional Lebesgue measure of its boundary is zero;
\item
\label{item:no-necks-dragging}
for every \(r\le \inr(\Omega)\) and every pair of balls \(B_r(x_0),B_r(x_1)\subset\Omega\), there exists a continuous curve \(\gamma:[0,1]\to\Omega\) such that
\[
\gamma(0)=x_0,
\qquad 
\gamma(1)=x_1,
\qquad 
B_r(\gamma(t))\subset\Omega,
\quad
\forall t\in[0,1].
\]
\end{enumerate}
\end{definition}

Condition~\ref{item:no-necks-regularity} imposes mild regularity on \(\partial\Omega\): the boundary must be a non-space-filling Jordan curve; we stress that there exist Jordan curves with positive two-dimensional Lebesgue measure (Osgood--Knopp curves; see~\cite{Sag94book}). Condition~\ref{item:no-necks-dragging} is a connectedness property at every scale \(r\): it says that one can continuously move a ball of radius \(r\) from one position to another inside \(\Omega\), without ever leaving the set.

For completeness, we note that sometimes only condition~\ref{item:no-necks-dragging} is taken as the definition of a set with no-necks, and condition~\ref{item:no-necks-regularity} is then imposed separately when one states the analogue of \cref{thm:LS22}.

The monotonicity of the isoperimetric profile with respect to the regularization parameter \(t\) in \cref{thm:monotonicity} implies the corresponding monotonicity of the Cheeger constant; see \cref{prop:fto}. In~\cite[Prop.~4.1]{Fto25}, the second author showed that this monotonicity may fail for non-convex sets. The example there is a set with no-necks, so \cref{thm:monotonicity} does not extend to this larger class. For completeness, we recall the relevant example.

\begin{example}[L-shaped set]
\label{ex:L-shaped}
Consider the set
\[
\Omega=(0,1)^2\cup [1,2)\times(0,\varepsilon),
\]
with \(\varepsilon>0\) sufficiently small. This is the reflection of the set considered in~\cite[Prop.~4.1]{Fto25} across the horizontal line \(y=\tfrac12\). Visually, it is an L-shaped set whose foot is very thin compared to its stem. For small values of the regularization parameter \(t\), the monotonicity of \(h(\widetilde{\Omega}_t)\) fails.

For completeness, we point out a minor computational mistake appearing in~\cite[Prop.~4.1]{Fto25}. The correct first-order expansion should read
\[
{|\Omega_{-t}|}\,h^2(\Omega_{-t})-{|\Omega|}\,h^2(\Omega)
=
(2+\sqrt{\pi})^2(4\varepsilon-2)t+O(t^2).
\]
In our volume-preserving notation, this becomes
\[
h^2(\widetilde{\Omega}_{-t})-h^2(\Omega)
=
\frac{1}{|\Omega|}
\left((2+\sqrt{\pi})^2(4\varepsilon-2)t+O(t^2)\right).
\]
\end{example}

In general, monotonicity may already fail at the level of the ball-like branch; in other words, Matheron’s inequality~\eqref{eq:matheron} may fail, that is, \(t\mapsto \inr(\widetilde{\Omega}_t)\) may not be monotonic. Let us look at its expansion near \(t=0^+\) for sets \(\Omega\) such that \(\inr(\Omega_t)=\inr(\Omega)+t\) for small \(t>0\), as is the case for the
example below. Since
\[
\inr(\widetilde{\Omega}_t)=\delta_\Omega(t)\bigl(\inr(\Omega)+t\bigr),
\]
the first-order Taylor expansion of \(\delta_\Omega(t)\) gives
\begin{equation}
\label{eq:matheron-failure}
\inr(\widetilde{\Omega}_t)
\ge
\inr(\Omega)
+
t\left(1-\frac12\frac{\Per(\Omega)}{|\Omega|}\inr(\Omega)\right)
+o(t),
\end{equation}
provided Steiner’s area formula~\eqref{eq:steiner-area} applies. The coefficient of the first-order term is non-negative for convex bodies by standard inequalities; see~\cite[Tab.~3.1, triplet \((A,p,r)\)]{SA00}, with equality for tangential bodies. For sets with no-necks it may be negative as the next example shows.

\begin{example}[A Pinocchio set]
\label{ex:pinocchio}
Consider the Pinocchio set introduced in~\cite[Ex.~4.6]{LP16}: it is the union of the unit ball centered at the origin, the ``face'', with a ball of radius \(\sin\theta\) centered at \((\cos\theta,0)\), the ``nose'', where \(\theta\in(0,\pi/2)\). Stretching the nose by a length \(\ell\), as in~\cite[Fig.~10 and Eq.~(4.2)]{LP16}, produces a set \(\mathcal P(\theta,\ell)\) whose perimeter and area are
\begin{align*}
\Per(\mathcal P(\theta,\ell))
&=
2(\pi-\theta)+\pi\sin\theta+2\ell,
\\
|\mathcal P(\theta,\ell)|
&=
(\pi-\theta)+\sin\theta\cos\theta+\frac{\pi\sin^2\theta}{2}+2\ell\sin\theta.
\end{align*}
Since \(\inr(\mathcal P(\theta,\ell))=1\), the first-order term coefficient in~\eqref{eq:matheron-failure} is negative if and only if
\[
\pi\sin\theta-2\sin\theta\cos\theta-\pi\sin^2\theta+2\ell(1-2\sin\theta) 
>
0.
\]
Choosing \(\theta<\pi/6\), one can find \(\ell=\ell(\theta)\gg 1\) such that this inequality holds. Geometrically, as the set is regularized, the nose shortens and thickens, and the face must initially shrink to compensate.
\end{example}

A further issue is that the regularized set \(\widetilde{\Omega}_t\) may fail to satisfy condition~\ref{item:no-necks-regularity}. The boundary may indeed cease to be parametrized by an injective curve. We give two examples, one for negative \(t\) and one for positive \(t\).

\begin{example}
\label{ex:breaking-pinocchio}
Consider again the Pinocchio set \(\mathcal P(\theta,\ell)\) of \cref{ex:pinocchio}, with \(\ell\gg0\). Choosing \(t=-\sin\theta\) produces an inner parallel set whose boundary consists of that of a Jordan domain together with a segment, which is crossed twice, once rightward and once leftward.
\end{example}

\begin{example}[U-shaped set]
\label{ex:U-shaped}
Consider the set
\[
\Omega
=
(0,1)^2
\cup
[1,2]\times(0,\varepsilon)
\cup
(2,3)\times(0,1).
\]
This is the union of the L-shaped set in \cref{ex:L-shaped} with its reflection across the vertical line \(x=\tfrac32\). For \(t=\tfrac12\), the outer parallel set \(\Omega_t\) has a non-injective boundary curve, since the vertical segment \(\{3/2\}\times[1/2+\varepsilon,1]\) is traversed once downward and once upward.
\end{example}

One possible workaround is to include an additional step in the regularization process: after performing the Minkowski sum or difference with a ball and rescaling, one could identify the set with the interior of its closure. This would  eliminate self-touching portions of the boundary and preserve the Jordan domain structure. However, continuity in \(t\) would then likely fail, and with it the joint continuity in \((v,t)\) from \cref{thm:reg-flow}.

A more serious issue, beyond self-touching, is a change of topology: the boundary may split into several Jordan curves, making it impossible to retain condition~\ref{item:no-necks-regularity}.

\begin{example}[C-shaped set]
Let \(C(\theta)\) be the closed convex cone
\[
C(\theta)=\{\,tp:\ t\ge 0,\ p\in \{\cos\theta\}\times[-\sin\theta,\sin\theta]\,\},
\]
and define
\[
\Omega=\bigl(B_{1}\setminus \overline{B}_{1-\varepsilon}\bigr)\setminus C(\theta).
\]
This set resembles the letter C, with two ends that are very close together if \(\theta\) is sufficiently small. Then, the outer parallel set \(\Omega_{\sin\theta}\) is not simply connected, hence it is not a Jordan domain. Indeed, it resembles an annulus, and its boundary is the union of two disjoint Jordan curves.
\end{example}

Summing up, for sets with no-necks, parallel bodies need not preserve the no-neck structure of \cref{def:no-neck}, and therefore one cannot directly exploit the analogue of \cref{thm:LS22} for all \(t\). We nonetheless expect the reparametrization law~\eqref{eq:shifted-RL} to survive in the range
\[
-\rho_\Omega<t<\reach(\Omega),
\]
where \(\reach(\Omega)\ge 0\) is the supremum of all non-negative \(s\) such that every point in \(\Omega_s\) has a unique projection onto \(\Omega\), see~\cite{Fed59,Tha08}, and \(\rho_\Omega\) denotes the largest \(r\) such that \(\Omega=\Omega_{-r}\oplus rB\). Indeed, for \(-\rho_\Omega<t<0\) one has \(\Omega=\Omega_t\oplus|t|B\), so that balls in \(\Omega\) of radius larger than \(|t|\) correspond to balls in \(\Omega_t\) as in \cref{lem:ball-to-ball}, and the dragging condition~\ref{item:no-necks-dragging} passes to \(\Omega_t\). For \(0<t<\reach(\Omega)\), every point of \(\Omega_t\setminus\Omega\) has a unique projection onto \(\Omega\); this yields \((\Omega\oplus tB)\ominus tB=\Omega\), hence the ball-to-ball correspondence, together with Steiner's formula for \(\Omega\), see~\cite{Fed59}.

The lower threshold cannot be improved: for the Pinocchio set of \cref{ex:pinocchio} one has \(\rho_{\mathcal P(\theta,\ell)}=\sin\theta\), and the no-neck structure breaks down precisely at \(t=-\sin\theta\), see \cref{ex:breaking-pinocchio}.


\bibliographystyle{plainurl} 
\bibliography{biblio.bib} 

\end{document}